\RequirePackage{fix-cm}
\documentclass[a4paper, 12pt, leqno]{article}%{amsart}
\pdfoutput=1

\usepackage{latexsym, amsmath, amsthm, amssymb, graphicx, graphics, epsfig, bm}
\usepackage{bm}
\usepackage{tikz-cd}

\usepackage{color}

\tikzset{cong/.style={draw=none,edge node={node [sloped, allow upside down, auto=false]{$\cong$}}},
         Isom/.style={every to/.append style={edge node={node [above,sloped, inner sep=0.4pt, allow upside down, auto=false]{$\sim$}}}}}

\usepackage{fullpage}

\usepackage{enumitem}

\newcommand{\comment}[1]{}

\usepackage[T1]{fontenc}
\usepackage{latexsym}
\usepackage{manfnt}
\usepackage[sans]{dsfont}
\usepackage{palatino}
\usepackage[sc, osf]{mathpazo} 
\usepackage{eulervm}
\usepackage{euscript}
\newlength{\otstup} 
\newcommand{\sub}[1]{\vspace{\otstup}\textbf{#1}\hspace*{0.5em}}
\newcommand{\upskip}{\vspace{-\otstup}}

\newtheorem{Thm}{Theorem}

\newtheorem{thm}{Theorem}[section]
\newtheorem*{thm*}{Theorem}
\newtheorem*{thmA'}{Theorem A'}
\newtheorem*{thmC'}{Theorem \ref{thm:K1}'}
\newtheorem*{thmD'}{Theorem \ref{thm:inv}'}
\newtheorem{lem}[thm]{Lemma}
\newtheorem*{lem*}{Lemma}

\newtheorem*{cor*}{Corollary}
\newtheorem{prop}[thm]{Proposition}
\newtheorem*{prop*}{Proposition}
\newtheorem*{claim*}{Claim}

\theoremstyle{definition}

\newtheorem{rem}[thm]{Remark}
\newtheorem*{rem*}{Remark}

\numberwithin{equation}{section}

\renewcommand{\proof}{\vspace{-6pt}\noindent\textit{\textbf{Proof. }}}

\renewcommand{\endproof}{$\square$\vspace{3pt}}

\newcommand{\txt}[1]{\;\;\text{#1}\;\;}
\newcommand{\q}[1]{{``#1''}}
\newcommand{\br}[1]{\left(#1\right)}
\newcommand{\brr}[1]{\left[#1\right]}
\newcommand{\bra}[1]{\left\langle #1 \right\rangle}

\newcommand{\isor}{\stackrel{\approx\hspace{0.2em}}{\longrightarrow}}

\newcommand{\hookto}{\hookrightarrow}
\newcommand{\hookr}{\arrow[hookrightarrow]}

\newcommand{\sett}[2]{\left\{#1 \;|\; #2\right\}}

\newcommand{\abs}[1]{\left| #1 \right|}

\newcommand{\inv}{^{-1}}
\newcommand{\restr}[2]{\left. #1 \right|_{#2}}

\renewcommand{\phi}{\varphi}

\newcommand{\matr}[1]{
\begin{pmatrix}
#1
\end{pmatrix}}

\newcommand{\smatr}[1]{\br{ \begin{smallmatrix}#1\end{smallmatrix} } }

\newcommand{\R}{\mathbb R}

\newcommand{\Z}{\mathbb Z}
\newcommand{\C}{\mathbb C}

\DeclareMathOperator{\ran}{ran}
\DeclareMathOperator{\coker}{coker}
\DeclareMathOperator{\Gr}{Gr}
\newcommand{\Grinf}{\Gr^{\star}}
\DeclareMathOperator{\LGr}{LGr}

\DeclareMathOperator{\Sympl}{Sympl}

\DeclareMathOperator{\codim}{codim}

\DeclareMathOperator{\GL}{GL}
\DeclareMathOperator{\Id}{Id}

\DeclareMathOperator{\dom}{{dom}}

\newcommand{\gr}{\Gamma}

\newcommand{\om}{^*}%{^{\omega}}
\newcommand{\Gom}{\Gamma\om}
\newcommand{\gammi}{\gamma^{-1}}

\newcommand{\sa}{^{\mathrm{sa}}}

\newcommand{\ort}{^{\mathrm{ort}}}

\comment{===============================================

=====================================================}

\newcommand{\GrF}{\Gr_{F}}

\newcommand{\U}{\mathcal U}

\newcommand{\UK}{\U_K}
\newcommand{\UF}{\U_F}

\newcommand{\Reg}{\mathcal R}
\newcommand{\Rsa}{\Reg\sa}
\newcommand{\RF}{\Reg_F}

\comment{===============================================
=====================================================}

\newcommand{\BFp}{\B_F^+}
\newcommand{\BFm}{\B_F^-}

\newcommand{\D}{\mathcal D}

\comment{===============================================

=====================================================}

\newcommand{\B}{\mathcal B}

\newcommand{\VV}{\mathcal V}

\newcommand{\BF}{\B_F}

\renewcommand{\H}{\mathcal H}

\comment{===============================================

=====================================================}

\newcommand{\xX}{x\in X}

\newcommand{\M}{\mathcal M}

\newcommand{\T}{\mathcal T}

\newcommand{\A}{\mathcal A}

\newcommand{\Hhat}{\hat{H}}
\newcommand{\Hh}{\Hhat}

\DeclareMathOperator{\ind}{\mathsf{ind}}

\newcommand{\kap}{\bm{\kappa}}
\newcommand{\kab}{\bar{\kap}}

\newcommand{\Dt}{\tilde{\mathcal D}}
\newcommand{\At}{\tilde{A}}
\newcommand{\AAt}{\tilde{\A}}
\newcommand{\Lt}{\widetilde{L}}
\newcommand{\mi}{^{\mathrm{min}}}
\newcommand{\ma}{^{\mathrm{max}}}
\newcommand{\Ami}{A}%{\A\mi}
\newcommand{\Ama}{A'}%{\A\ma}
\newcommand{\Op}{A}

\comment{===============================================

=====================================================}

\newcommand{\dG}{\partial G}

\newcommand{\CDS}{\mathcal{C}}
\newcommand{\cds}{C}

\title{Family index for Fredholm extensions \\ of semi-Fredholm operators} 
\author{Marina\,Prokhorova}
\date{}

\begin{document}

\maketitle

\footnotetext{\hspace*{-1.8em}Department of Mathematics, University of Haifa (Israel) \\
Department of Mathematics, Technion -- Israel Institute of Technology \\
ORCID 0000-0003-3214-719X}

\renewcommand{\baselinestretch}{1.05}
\selectfont

\upskip
\begin{abstract}
\parindent=0cm
\setlength{\parskip}{0pt plus 0pt minus 0pt}
\noindent
This paper is devoted to an abstract analogue of elliptic boundary value problems,
namely, Fredholm realizations of semi-Fredholm operators in a Hilbert space.
Such a realization is determined by an abstract boundary condition,
which is a subspace in the space of abstract boundary values.
We find the $K^0$ index of a family of such abstract boundary value problems,
or the $K^1$ index in the self-adjoint case,
in terms of the corresponding family of abstract boundary conditions.
Our approach is based on passing from a Fredholm operator to its graph.
The graph forms a Fredholm pair with the horizontal subspace,
and we prove the index formula by deforming the horizontal subspace instead of the operator.
\end{abstract}

\let\oldnumberline=\numberline
\renewcommand{\numberline}{\vspace{-19pt}\oldnumberline}
\addtocontents{toc}{\protect\renewcommand{\bfseries}{}}

%\upskip
{\small\tableofcontents}
\bigskip

\section{Introduction}\label{sec:intro}

\upskip
\sub{Unbounded operators.}
Let $H$ be a separable complex Hilbert space of infinite dimension
and $A$ be a closed densely defined linear operator on $H$.
The graph $\gr_A$ of $A$ is a closed subspace of $H\oplus H$ and thus is itself a Hilbert space.
It can be identified with the domain $\dom A$ of $A$ equipped with the graph inner product
\[ \bra{\xi,\eta}_A = \bra{\xi,\eta} + \bra{A\xi,A\eta} \txt{for} \xi,\eta\in\dom A. \]
We will always consider domains of unbounded operators with this Hilbert space structure.

The graph $\Gamma_A$ of a closed operator $A$ is a closed subspace of $H\oplus H$.
For our purposes, an appropriate topology on the set of closed operators 
is the \emph{graph topology}.
By the definition, a family of operators is continuous in the graph topology
if the family of orthogonal projections onto their graphs is norm continuous.
For self-adjoint operators, this is the same as the uniform resolvent continuity.

A graph continuous family $A=(A_x)_{\xX}$ of Fredholm operators has the homotopy invariant 
taking values in the K-group $K^0(X)$ and called the index of $A$.
If all operators $A_x$ are self-adjoint, then the relevant invariant takes values in $K^1(X)$.

A continuous family of linear elliptic operators of positive order on a smooth manifold with boundary,
taken together with a continuous family of elliptic boundary conditions,
leads to a graph continuous family of unbounded Fredholm operators.
This paper is devoted to an abstract analogue of this situation 
and to the family index for abstract boundary value problems.

\sub{Extension pairs and their realizations.} 
A closed operator $\Ama$ is called an \textit{extension} of $\Ami$, written as $\Ami\subset\Ama$, 
if $\dom\Ami\subset\dom\Ama$ and the restriction of $\Ama$ to $\dom\Ami$ coincides with $\Ami$.
We call such a pair $(\Ami,\Ama)$ an \textit{extension pair}.
A \textit{realization} of an extension pair $\A=(\Ami,\Ama)$ is a closed operator $\At$ 
such that $\Ami\subset\At\subset \Ama$.
It is given by the restriction of $\Ama$ to the domain $\Dt=\dom\At$, 
which is a closed subspace of $\D'=\dom\Ama$ containing $\D=\dom\Ami$.
Therefore, realizations of $\A$ are in one-to-one correspondence with such subspaces $\Dt$.

Passing from $\Dt$ to the quotient $L = \Dt/\D \subset \D'/\D$
provides a one-to-one correspondence between realizations of $\A$ 
and closed subspaces $L$ of the quotient space $\VV = \D'/\D$.
This quotient space is called the \textit{space of (abstract) boundary values} of $\A$,
and the quotient map 
\[ \gamma\colon\dom\Ama\to\dom\Ama/\dom\Ami=\VV \] 
is called the \textit{(abstract) trace map}.
A closed subspace $L$ of $\VV$ is called an \textit{(abstract) boundary condition} for $\A$;
it corresponds to the realization $\A_L$ with the domain 
\[ \dom\A_L = \sett{\xi\in\dom\Ama}{\gamma\xi\in L}. \]

\upskip
\sub{Fredholm extension pairs.}
Suppose that an extension pair $\A=(\Ami, \Ama)$ admits a Fredholm realization $\At$.
Then the restriction $\Ami$ of $\At$ is lower semi-Fredholm
(that is, its range is closed and its kernel is finite-dimensional)
and the extension $\Ama$ of $\At$ is upper semi-Fredholm
(that is, its range is closed of finite codimension).
We call such a pair $(\Ami, \Ama)$ a \textit{Fredholm extension pair}.
Conversely, let $\A=(\Ami, \Ama)$ be a Fredholm extension pair.
Then it admits Fredholm realizations. 

Fredholm extension pairs and their realizations appear naturally in applications.
For example, if $M$ is a compact manifold with boundary
and $D$ is an elliptic operator of positive order on $M$, 
then its minimal realization $D\mi$ and maximal realization $D\ma$ form a Fredholm extension pair 
$D\mi\subset D\ma$.
Realizations of this extension pair are given by boundary conditions for $D$.
Elliptic boundary conditions give rise to Fredholm realizations.

\sub{Cauchy data space.}
Let $(\Ami, \Ama)$ be a Fredholm extension pair.
In order to determine, in terms of a boundary condition $L$, 
whether its realization $\A_L$ is Fredholm,
one should compare $L$ with a distinguished closed subspace $\CDS$ of $\VV$. 
It is defined as the image 
\[ \CDS=\gamma(\ker A') \] 
of the kernel of $\Ama$ under the trace map $\gamma\colon\dom(\Ama)\to\VV$
and is called the \textit{(abstract) Cauchy data space} of the extension pair $\A$.
The corresponding realization $\A_\CDS$ has the domain $\dom\Ami+\ker\Ama$ and 
is characterized by the properties $\ker\A_\CDS=\ker\Ama$ and $\ran\A_\CDS=\ran\Ami$,
so it has the maximal possible kernel and the minimal possible range.

Recall that a pair $(L,M)$ of closed subspaces of a Hilbert space $\VV$ is called \textit{Fredholm} 
if $L+M$ is closed and both $\dim(L\cap M)$ and $\codim(L+M)$ are finite.
Such a pair has an integer-valued invariant called the index,
\[ \ind_{\VV}(L,M) = \dim(L\cap M)-\codim(L+M). \] 
(The subscript $\VV$ at $\ind$ points out to the space where the Fredholm pair $(L,M)$ lives.)

Let $\A$ be a Fredholm extension pair and $L$ be a closed subspace of $\VV$.
Then the realization $\A_L$ is Fredholm if and only if the pair $(L,\CDS)$ is Fredholm,
and in this case the index of the realization $\A_L$ is given by the formula
	\begin{equation}\label{eq:FGS}
		\ind\A_L = \ind_{\VV}(L,\CDS) + \dim\ker\Ami - \dim\coker\Ama \in\Z. 
	\end{equation}
See \cite[Theorem B.28]{BGS}.
The first part of this paper is devoted to a generalization of this formula to the family case.

\sub{Families of boundary conditions.}
Suppose that a boundary condition $L_x\subset\VV$ continuously depends on a parameter $\xX$,
where $X$ is a topological space.
(Here \q{continuity} means that the orthogonal projection onto $L_x$ depends continuously on $x$
in the norm topology.) 
Then the corresponding family of realizations $\AAt_x=\A_{L_x}$ is continuous in the graph topology 
(that is, the graph of $\AAt_x$ depends continuously on $x$).
See Proposition \ref{prop:Gr-Reg-cont}.

Suppose additionally that each realization $\AAt_x$ is Fredholm.
As was mentioned above, a graph continuous family of Fredholm operators has a homotopy invariant 
$\ind\AAt\in K^0(X)$. 
A continuous family $(L_x,\CDS)$ of Fredholm pairs in $\VV$ also
has a homotopy invariant $\ind_{\VV}(L,\CDS)\in K^0(X)$. % called the index. 
See Section \ref{sec:K0ind} for details.

The following theorem presents the simplest case of our first result.

\begin{Thm}\label{Thm:ind0}
Let $(L_x)$ be a continuous family of Fredholm boundary conditions 
for a Fredholm extension pair $\A = (\Ami,\Ama)$.
Suppose that $\Ami$ is injective and $\Ama$ is surjective. Then 
\begin{equation*}%\label{eq:Thm-ind-K0}
	\ind\A_L = \ind_{\VV}(L,\CDS) \in K^0(X). 
\end{equation*}
\end{Thm}

\upskip
\sub{Families of extension pairs.}
More generally, an extension pair may also depend on a parameter.
An appropriate notion of continuity in this context is graph continuity,
both of the minimal operator $A$ and maximal operator $A'$.

Let $\A_x = (\Ami_x,\Ama_x)$, $\xX$ be a \emph{graph continuous} family of Fredholm extension pairs
parametrized by points of $X$.
It gives rise to a continuous family of spaces of boundary values $\VV_x = \dom\Ama_x/\dom\Ami_x$.
Let $L_x\subset\VV_x$ be a Fredholm boundary condition for $\A_x$ depending continuously on $x$. 
Then the corresponding family of Fredholm realizations $\A_L=(\A_{x,L_x})_{\xX}$ is graph continuous,
so its index $\ind\A_L\in K^0(X)$ is well defined.

The index formula \eqref{eq:FGS} for a single operator includes the dimensions 
of the kernel of $\Ami$ and the cokernel of $\Ama$.
If $A$ and $A'$ depend on a parameter, then the corresponding families of subspaces 
are not necessarily continuous and might not determine $K^0$ classes. 
Because of this, in the present paper we consider only a special case,
namely, we suppose that the dimensions of $\ker \Ami_x$ and $\coker \Ama_x$ are \textit{locally constant}.
Then the finite-dimensional spaces $\ker \Ami_x$ and $\coker \Ama_x$ depend continuously on $x$, 
so they represent elements $[\ker\Ami]$ and $[\coker\Ama]$ of $K^0(X)$.
Moreover, in this case the Cauchy data space $\CDS_x$ of $\A_x$ also depends continuously on $x$. 

\begin{Thm}\label{Thm:ind-AL}
Let $\A_x = (\Ami_x,\Ama_x)$, $\xX$ be a graph continuous family of Fredholm extension pairs
such that $\dim\ker\Ami_x$ and $\dim\coker\Ama_x$ are locally constant.
Let $L_x$ be a continuous family of Fredholm boundary conditions for $\A_x$.
Then 
\begin{equation}\label{eq:Thm-ind-K0}
	\ind\A_L = \ind_{\VV}(L,\CDS) + [\ker\Ami] - [\coker\Ama] \in K^0(X). 
\end{equation}
In particular, $\ind\A_L = \ind_{\VV}(L,\CDS)$ if each $\Ami_x$ is injective and each $\Ama_x$ is surjective.
\end{Thm}

\upskip
We prove Theorems \ref{Thm:ind0} and \ref{Thm:ind-AL} in Section \ref{sec:ind-ext}.

If the space $X$ consists of one point, then $\A_L$ is a single operator, its index is integer-valued, 
and equality \eqref{eq:Thm-ind-K0} takes the form \eqref{eq:FGS}.

\sub{Self-adjoint extensions.} 
The classical von Neumann theory of self-adjoint extensions of a symmetric operator 
may be formulated in two equivalent ways:
in terms of unitary operators or in terms of symplectic geometry and Lagrangian subspaces. 
See \cite[Section XII.4]{DS} or \cite[Section X.1]{RS}.
We will use symplectic language as it is more appropriate for our purposes.

Let $A$ be a closed densely defined symmetric operator on $H$ and $A^*$ be the operator adjoint to $A$.
Every symmetric extension of $A$ is a restriction of $A^*$, 
so symmetric extensions of $A$ are realizations of the extension pair $\A=(A,A^*)$.

The space $\VV_A = \dom A^*/\dom A$ of abstract boundary values of a symmetric operator $A$ 
is equipped with a sesquilinear form 
\begin{equation}\label{eq:omega-A}
	\omega(\gamma\xi,\gamma\eta) = \bra{A^*\xi,\eta} - \bra{\xi,A^*\eta} 
		\quad\text{for } \xi,\eta\in\dom A^*, 
\end{equation}
where $\gamma\colon\dom A^* \to \VV$ is the trace map.
This form turns $\VV$ into a (complex) symplectic space.
For a subspace $L\subset\VV$, the extension $A_L=\A_L$ of $A$ is the restriction of $A^*$ 
to the domain $\dom A_L = \gamma\inv(L)$.
It is symmetric if and only if $L$ is isotropic,
and self-adjoint if and only if $L$ is Lagrangian.
Therefore, self-adjoint extensions of $A$ are parametrized by Lagrangian subspaces of $\VV$.

\sub{Fredholm self-adjoint extensions.}
A symmetric operator $A$ admits Fredholm self-adjoint extensions if and only if 
the  extension pair $(A,A^*)$ is Fredholm or, equivalently, the operator $A$ is semi-Fredholm
(we omit the word \q{lower} since a semi-Fredholm symmetric operator is necessarily lower semi-Fredholm).

Let $A$ be a semi-Fredholm symmetric operator.
Then the Cauchy data space $\CDS = \gamma(\ker A^*)$  is a Lagrangian subspace of $\VV$
and Fredholm self-adjoint extensions of $A$ are parametrized 
by points of the Fredholm Lagrangian Grassmanian 
\[ \LGr_F(\VV;\CDS) = \sett{L\subset\VV}{L \text{ is Lagrangian and } (L,\CDS) \text{ is a Fredholm pair}}. \]

Suppose that we have a continuous family $L_x$ of Fredholm self-adjoint boundary conditions for $A$ 
parametrized by points of a topological space $X$.
It gives rise to a graph continuous family of Fredholm self-adjoint extensions $A_L = (A_{L_x})_{\xX}$
(equivalently, this family is continuous in the uniform resolvent sense). 
The family $A_L$ has a homotopy invariant called the index which takes values in $K^1(X)$.
The family $(L_x,\CDS)$ of Fredholm Lagrangian pairs in $\VV$ also has a homotopy invariant 
$\ind_{\VV}(L,\CDS)$ taking values in $K^1(X)$.
See Section \ref{sec:K1ind} for details.

\begin{Thm}\label{Thm:K1-AL}
Let $L_x$ be a continuous family of Fredholm self-adjoint boundary conditions 
for a semi-Fredholm symmetric operator $A$. 
Then $\ind A_L = \ind_{\VV}(L,\CDS)\in K^1(X)$.
\end{Thm}

More generally, an operator $A$ may also depend on a parameter $\xX$.
We suppose that the dimension of the kernel of $A_x$ is locally constant;
then the Cauchy data space $\CDS_x$ of $A_x$ depends continuously on $x$.

\begin{Thm}\label{Thm:K1-AL-k}
Let $A_x$, $\xX$ be a graph continuous family of semi-Fredholm symmetric operators 
such that $\dim\ker A_x$ is locally constant.
Let $L_x$ be a continuous family of Fredholm self-adjoint boundary conditions for $A_x$.
Then $\ind A_L = \ind_{\VV}(L,\CDS)\in K^1(X)$.
\end{Thm}

\upskip
We prove Theorems \ref{Thm:K1-AL} and \ref{Thm:K1-AL-k} in Section \ref{sec:K1-op};
see also Theorems \ref{thm:K1-V} and \ref{thm:K1-Vx}.

A particular case of these results is the \q{spectral flow = Maslov index} formula for loops.
Indeed, if $X$ is a circle, then, under the identification $K^1(S^1)\cong\Z$, 
the index of the loop of Fredholm self-adjoint operators coincides with its spectral flow,
while the index of the loop of Fredholm Lagrangian pairs is its Maslov index. 

Theorem \ref{Thm:K1-AL-k} is related to a theorem from a recent work of N. Ivanov.
See \cite{I23}, Theorem 5.1.
After learning about Theorem \ref{Thm:K1-AL-k}, Ivanov deduced in \cite{I23b} 
an essentially equivalent result from his Theorem 5.1.

It is sometimes more convenient to state self-adjoint boundary conditions in terms of boundary triplets.
In section \ref{sec:triplets} we translate Theorems \ref{Thm:K1-AL} and \ref{Thm:K1-AL-k}
to this language, see Theorems \ref{thm:K1-triplet} and \ref{thm:K1-triplet-inv}.
We illustrate Theorem \ref{thm:K1-triplet} by applying it to 
Schr\"odinger operator on finite metric graphs.
See Theorem \ref{thm:index-QG}.
A particular case of Theorem \ref{thm:index-QG}, 
namely the \q{spectral flow = Maslov index} formula for loops of boundary conditions on such a graph, 
was recently proven in \cite{BPS} by different methods. 
Our result generalizes \cite[Theorem A]{BPS} to arbitrary families of boundary conditions.

\sub{Elliptic operators.}
Theorems \ref{Thm:ind0}--\ref{Thm:K1-AL-k} have natural applications to elliptic boundary value problems.
A continuous family of elliptic differential operators on a compact manifold with boundary 
leads to graph continuous families of their minimal and maximal realizations,
so it falls within the framework of Theorems \ref{Thm:ind-AL} and \ref{Thm:K1-AL-k}.

The space of abstract boundary values of such an operator $D$ 
is naturally identified with a mix of Sobolev spaces, both positive and negative fractional orders, 
over the boundary, where the mix depends on the boundary operator of $D$.
It would be more convenient to deal with boundary conditions given by subspaces in a single Sobolev space, for example in the $L^2$-space of sections over the boundary.
We provide such an adaptation of Theorems \ref{Thm:ind0}--\ref{Thm:K1-AL-k} 
in the next paper.

\sub{Passing to Fredholm pairs.}
In our proof of Theorems \ref{Thm:ind0}--\ref{Thm:K1-AL-k} we first pass from operators to their graphs.
A Fredholm operator $\At$ is replaced by its graph $\gr_{\At}\subset \Hh=H\oplus H$,
which is Fredholm with respect to the horizontal subspace $H\oplus 0$.
This provides an embedding (that is, a homeomorphism on the image)
of the space of Fredholm operators equipped with the graph topology to the Fredholm Grassmanian 
\[ \Gr_F(H\oplus H; H\oplus 0) = \sett{\Gamma\subset H\oplus H}{(\Gamma,H\oplus 0) \text{ is a Fredholm pair}}. \]
By a previous result of the author, 
this embedding is a homotopy equivalence.
See \cite[Theorem A]{Pr21}.
In particular, the $K^0$-index of the family $\At_x$ is equal 
to the homotopy class of the family $\gr_{\At_x}$ in $\Gr_F$.
Equivalently, it is equal to the index of the corresponding family of Fredholm pairs $(\gr_{\At_x},H\oplus 0)$:
\begin{equation}\label{eq:ind=}
	\ind\At = \ind_{H\oplus H}(\gr_{\At},H\oplus 0) \in K^0(X). 
\end{equation}
In such a way, the study of Fredholm realizations of a Fredholm extension pair $(\Ami,\Ama)$ 
is replaced by the study of subspaces of $H\oplus H$ which contain the graph $\Gamma$ of $\Ami$, 
are contained in the graph $\Gamma'$ of $\Ama$,
and are Fredholm with respect to the horizontal subspace $H\oplus 0$.
This provides more possibilities for homotopies than working with operators only.
Moreover, the index of a family of Fredholm pairs $(\gr_{\At},H\oplus 0)$ is invariant
under continuous deformations not only of the first argument $\gr_{\At}$, 
but also of the second argument $H\oplus 0$.
We use this property in the second part of our proof. 

A similar reasoning works for self-adjoint operators. 
In this case, the Fredholm Grassmanian should be replaced by the Fredholm Lagrangian Grassmanian,
$K^0$ by $K^1$, and Theorem A from \cite{Pr21} by Theorem C.

\sub{Index formula for Fredholm pairs.}
The second part of our proof is provided by an index formula for Fredholm pairs of certain form.
See Theorem \ref{thm:ind-Gr} for the non-self-adjoint case and 
Theorem \ref{thm:ind-LGr} for the self-adjoint case.
These results have an independent interest, so we present Theorem \ref{Thm:ind-Gr1}, 
which is a particular case of Theorem \ref{thm:ind-Gr}, here in the Introduction. 
Its $K^1$ analogue, Theorem \ref{thm:ind-LGr1}, is proven in the same manner. 

Let $\Gamma\subset\Gamma'$ be closed subspaces of a Hilbert space $\Hhat$ and 
\[ \gamma\colon\Gamma'\to\Gamma'/\Gamma=\VV \] 
be the quotient map. 
The inverse image $\gamma\inv(L)$ of a closed subspace $L\subset\VV$ 
is a closed subspace of $\Gamma'\subset\Hhat$ containing $\Gamma$.

Let $M$ be a closed subspace of $\Hhat$ such that 
\[ M+\Gamma \text{ is closed}, \quad M\cap \Gamma = 0, \txt{and} M+\Gamma'=\Hhat \]
(in this case we say that $M$ is \textit{transversal} to the pair $(\Gamma,\Gamma')$).
Then the pair $(\gamma\inv(L),M)$ is Fredholm in $\Hhat$ if and only if 
the pair $(L,\gamma(M\cap \Gamma'))$ is Fredholm in $\VV$. 
See Proposition \ref{cor:gamma-LM-tF}.

\begin{Thm}\label{Thm:ind-Gr1}
Let $M$ be transversal to the pair $(\Gamma\subset\Gamma')$.
	Let $L=(L_x)_{\xX}$ be a continuous family of subspaces of $\VV$ 
	which are Fredholm with respect to $\CDS=\gamma(M\cap \Gamma')$.
	Then  
	\begin{equation}\label{eq:Thm:ind-Gr1}
		\ind_{\Hhat}(\gamma\inv(L),M) = \ind_{\VV}(L,\CDS) \in K^0(X). 
	\end{equation}
\end{Thm}

If $\Hh=H\oplus H$, $M=H\oplus 0$, $\Gamma$ is the graph of $\Ami$, and $\Gamma'$ is the graph of $\Ama$,
then $\gamma\inv(L)$ is the graph of $\A_L$ and $\gamma(M\cap \Gamma') = \gamma(\ker A')=\CDS$
(this explains our use of notation $\CDS$ for $\gamma(M\cap \Gamma')$).
Combining Theorem \ref{Thm:ind-Gr1} with \eqref{eq:ind=}, we obtain Theorem \ref{Thm:ind0}.

Our proof of Theorem \ref{Thm:ind-Gr1} is based on a continuous deformation 
$\M_t$ of $\M_0=M$ to $\M_1 = \CDS\oplus F$, where $F = \Hhat\ominus\Gamma'$,
such that $\gamma(\M_t\cap \Gamma') = \CDS$ for all $t\in[0,1]$.
The pair $(\gamma\inv(L_x),\M_t)$ is Fredholm for every value of parameters,
so the homotopy invariance of the index implies 
\[ \ind_{\Hhat}(\gamma\inv(L),M) = \ind_{\Hhat}(\gamma\inv(L),\CDS\oplus F). \]
Additivity of the index with respect to direct sums implies
\[ \ind_{\Hhat}(\gamma\inv(L),\CDS\oplus F) = 
   \ind_{\VV}(L,\CDS) + \ind_{\Gamma\oplus F}(\Gamma, F) = \ind_{\VV}(L,\CDS), \]
since the index of the pair $(\Gamma, F)$ in $\Gamma\oplus F$ vanishes.
Taking the last two formulas together, we obtain equality \eqref{eq:Thm:ind-Gr1}.

Our construction of the homotopy $\M_t$ provided by Proposition \ref{prop:TN} is canonical.
Therefore, it works in families as well, when not only $L$ 
but also all the other variables $\Gamma$, $\Gamma'$, $M$, $\Hh$ depend on a parameter $\xX$.
This provides Theorem \ref{thm:ind-Gr}, a parametrized version of Theorem \ref{Thm:ind-Gr1}, 
from which Theorem \ref{Thm:ind-AL} is deduced.

\sub{The case of compact base space.}
The reasoning above allows us to prove all the theorems for an arbitrary base space $X$.
In the particular case of a \textit{compact} space $X$,
Theorems \ref{Thm:ind0}, \ref{Thm:ind-AL}, and \ref{Thm:ind-Gr1} can be proved in a more standard way, 
using a deformation of the family of boundary conditions $L$ to a family $L'$ 
for which $L'\cap\CDS$ and $L'+\CDS$ depend continuously on a parameter.
See Section \ref{sec:alt} for details.

This standard approach, however, does not work for the $K^1$ index, even if $X$ is compact.
Indeed, a family $L$ admits such a deformation in the class of Fredholm self-adjoint boundary conditions 
only if $\ind(L,\CDS) = 0\in K^1(X)$.

\sub{Acknowledgments.} 
This work was partially supported by the ISF grants no. 431/20 and 2848/25 
and by the European Research Council (ERC) under the European Union's Horizon 2020 research and innovation program (grant no.~101001677).
It was partially done during my postdoctoral fellowships at the Technion -- Israel Institute of Technology
and at the University of Haifa and during my visit to the Max Planck Institute for Mathematics in Bonn.

I am grateful to N. V. Ivanov for reading a preliminary version of this paper 
and helping to improve the exposition.

\addtocontents{toc}{\vspace{\otstup}}

\section{Grassmanian}\label{sec:Gr}  

In this section we recall some (mostly standard) notions and facts that will be used throughout the paper.

\sub{Grassmanian.}
Let $H$ be a complex Hilbert space.
The set $\Gr(H)$ of closed subspaces of $H$ is called the \textit{Grassmanian}.
It is equipped with the topology induced from the norm topology on the space $\B(H)$ of bounded operators
by the inclusion $\Gr(H)\hookto\B(H)$ taking a closed subspace $M\subset H$ 
to the orthogonal projection $P_M$ onto $M$.

It is sometimes convenient to deal with arbitrary idempotents instead of orthogonal projections 
(that is, self-adjoint idempotents).

\begin{prop}\label{prop:idemp}
	If $T_x\in\B(H)$ is a family of idempotents depending continuously on $x$,
	then the ranges of $T_x$ form a continuous family of closed subspaces of $H$,
	and similarly for the kernels of $T_x$. 
\end{prop}

\proof
For every idempotent $T$, the operator $T+T^*-1$ is invertible and the orthogonal projection $T\ort$ 
onto the range of $T$ is given by the formula $T\ort = T(T+T^*-1)\inv$. 
See \cite[Proposition A.2(2)]{Pr17}.
Therefore, the map $T\mapsto T\ort$ is continuous 
and thus also the map taking $T$ to $\ran T = \ran T\ort$ is continuous
(here $\ran T$ denotes the range of $T$).

The kernel of $T_x$ is the range of the idempotent $1-T_x$ and thus also depends continuously on $x$.
\endproof

\begin{prop}\label{prop:Gr2Gr-b}
	A bijective operator $F\colon H\to H'$ induces a homeomorphism 
	\[ F\colon\Gr(H) \isor \Gr(H'), \quad M\mapsto F(M). \] 
\end{prop}

\proof
By the Closed Graph Theorem, both $F$ and $F\inv$ are bounded.
The subspace $F(M)$ is the range of the idempotent $FP_MF\inv$, which depends continuously on $M$.
It remains to apply Proposition \ref{prop:idemp}.
\endproof

\sub{Canonical decomposition.} 
Recall a standard construction that we will use throughout the paper.
The trivial Hilbert bundle over $\Gr(H)$ with the fiber $H$ is canonically decomposed into the direct sum 
	$H\times\Gr(H) = \H\oplus\H^{\bot}$
of two Hilbert bundles, whose fibers over $M\in\Gr(H)$ are $\H_M = M$ and $\H^{\bot}_M = M^{\bot}$,
the orthogonal complement of $M$ in $H$.
This decomposition is locally trivial in the following sense: for every $M_0\in\Gr(H)$ 
there is a continuous map $g\colon V\to\U(H)$ from a neighborhood $V$ of $M_0$ to the unitary group 
such that $M = g_M M_0$ for every $M\in V$.
See \cite[Proposition 5.2.6]{WO}.

An important subspace of $\Gr(H)$ is
\[ \Grinf(H) = \sett{M\in\Gr(H)}{\dim M = \infty = \codim M}. \]
The restrictions of $\H$ and $\H^{\bot}$ to $\Grinf(H)$ are locally trivial Hilbert bundles 
over a paracompact space, 
and their structure group is the unitary group $\U$ with the norm topology, which is contractible 
by Kuiper's theorem \cite{Kui}. 
Thus these restrictions are trivial as Hilbert bundles with the structure group $\U$.
This can be stated as follows. 

\begin{prop}\label{prop:decomp}
	For a subspace $M_0\in\Grinf(H)$, there is a map 
	\begin{equation}\label{eq:gP*U}
		g\colon \Grinf(H)\to\U(H) \;\text{ such that }\; M = g_M M_0 \;\text{ for every }\; M\in\Grinf(H),
	\end{equation}
\end{prop}

The map $\U(H)\to\Grinf(H)$, $u\mapsto uM_0$ has a section \eqref{eq:gP*U} 
and thus is a trivial bundle with contractible structure group $\U(M_0)\times\U(M_0^{\bot})$
and contractible total space $\U(H)$ over the  paracompact space $\Grinf(H)$.
It follows that $\Grinf(H)$ is also contractible.
See \cite[proof of Lemma 3.6]{ASi}.

\sub{Maps with changing domains.}
Let $M_x$ and $N_x$ be continuous families of closed subspaces of $H$.
Then the notion of continuity for a family of bounded operators $T_x\colon M_x\to N_x$ is well defined.
Indeed, locally one can chose a continuous family $u_x$ of isometries 
from some fixed Hilbert space $M$ to $H$ such that the range of $u_x$ is $M_x$.
Combining $T_xu_x\colon M\to N_x$ with the inclusion $N_x\hookto H$, 
one gets the family of operators $M\to H$ between fixed Hilbert spaces, 
where we have the standard notion of norm continuity.
Equivalently, one can replace the family $T_x$ with the family of operators $T_x\oplus 0\in\B(H)$,
where $0$ denotes the zero operator $M_x^{\bot}\to N_x^{\bot}$. 

\begin{rem}\label{rem:chang-dom}
	It can be easily seen that analogues of Propositions \ref{prop:idemp} and \ref{prop:Gr2Gr-b} 
	hold in this situation as well.
	Namely, a continuous family $T_x\in\B(M_x)$ of idempotents gives rise 
	to a continuous family of closed subspaces $\ran(T_x)$ and $\ker(T_x)$ of $M_x$ (and thus of $H$).
	A continuous family of closed subspaces $L_x\subset M_x$
	together with a continuous family of isomorphisms $F_x\colon M_x\to N_x$ 
	gives rise to a continuous family of closed subspaces $F_x(L_x)$. 
\end{rem}

\upskip
\sub{Sum and intersection of subspaces.}
The following result is proved in \cite{Neu} for Banach spaces; 
we need it only for Hilbert spaces.

\begin{prop}[\cite{Neu}, Lemma 1.5]
\label{prop:cap-cup-cont}
Let $M_x$ and $N_x$ be continuous families of closed subspaces of $H$ 
such that $M_x+N_x$ is closed for every $x$.
Then $M_x+N_x$ depends continuously on $x$ if and only if $M_x\cap N_x$ depends continuously on $x$.
\end{prop}

\proof
First note that if $M_x$ and $N_x$ are orthogonal, then the orthogonal projection onto 
$M_x+N_x = M_x\oplus N_x$ is equal to the sum of orthogonal projections onto $M_x$ and $N_x$,
so the sum $M_x+N_x$ depends continuously on $x$.
Similarly, if $C_x\subset M_x$ is a continuous family of subspaces, 
then the orthogonal complement $M_x\ominus C_x$ depends continuously on $x$.

By \cite[Lemma IV.4.9]{Kato}, the sum $M^{\bot}+N^{\bot}$ is closed if and only if $M+N$ is closed.
The map taking a subspace to its orthogonal complement is a homeomorphism of $\Gr(H)$; 
it replaces (closed) sums with intersections and vice versa.
Therefore, it is enough to prove that continuity of $C_x=M_x\cap N_x$ implies continuity of $M_x+N_x$.

Passing from the pair $(M_x,N_x)$ to the pair $(M'_x,N'_x) = (M_x\ominus C_x,N_x\ominus C_x)$, 
we get two continuous families $M'_x$ and $N'_x$ of subspaces with trivial intersections $M'_x\cap N'_x = 0$ 
and closed sums $L_x = M'_x+N'_x = (M_x+N_x)\ominus C_x$.
The continuity of $L_x$ would imply continuity of the orthogonal sum $L_x\oplus C_x = M_x+N_x$.
Thus we need only to show that $M_x\cap N_x \equiv 0$ implies continuity of $M_x+N_x$ provided these sums are closed.

Suppose that $M_x\cap N_x = 0$ for all $x$, and let $z\in X$. 
We will show that $M_x+N_x$ depends continuously on $x$ in some neighborhood of $z$. 
Choose a complement subspace $V\subset H$ for $M_z+N_z$.
The bounded operator $S_x\colon M_x\oplus N_x\oplus V\to H$ taking $m\oplus n\oplus v$ to $m+n+v$
depends continuously on $x$ and is invertible for $x=z$. 
Therefore, it is invertible also for $x$ in some neighborhood $U$ of $z$.
The subspace $M_x\oplus N_x$ of $M_x\oplus N_x\oplus V$ depends continuously on $x$,
so Proposition \ref{prop:Gr2Gr-b} together with Remark \ref{rem:chang-dom}
imply continuity of the family $S_x(M_x\oplus N_x) = M_x+N_x$ for $x\in U$. 
\endproof

\section{Regular operators and Fredholm pairs of subspaces}\label{sec:GrF}

\upskip
\sub{Fredholm and transversal pairs.}
The concepts of Fredholm and semi-Fredholm pairs of subspaces was introduced by Kato \cite{Kato}.
A pair $(M,N)$ of closed subspaces of a Hilbert space $\Hhat$ is called: 
\begin{itemize}[topsep=-2pt, itemsep=0pt, parsep=3pt, partopsep=0pt]
	\item \textit{transversal} if $M\cap N = 0$ and $M+N = \Hhat$;
	\item \textit{Fredholm} if $M+N$ is closed and $M\cap N$, $\Hhat/(M+N)$ are finite-dimensional;
	\item \textit{lower semi-Fredholm} if $M+N$ is closed and $M\cap N$ is finite-dimensional;
	\item \textit{upper semi-Fredholm} if $M+N$ is closed and finite-codimensional. 
\end{itemize}

A pair $(M,N)$ is lower semi-Fredholm if and only if $(M^{\bot},N^{\bot})$ is upper semi-Fredholm.
A pair $(M,N)$ is Fredholm if and only if $(M^{\bot},N^{\bot})$ is Fredholm.
See \cite[Corollary 4.13]{Kato}.
The same argument shows that $(M,N)$ is transversal if and only if $(M^{\bot},N^{\bot})$ is transversal.

We will use two kinds of the Fredholm Grassmanian:
\[ \Gr_F(\Hhat;M) = \sett{L\in\Gr(\Hhat)}{(L,M) \text{ is a Fredholm pair}}, \] 
\[ \Gr_F(\Hhat) = \sett{(L,M)\in\Gr(\Hhat)^2}{(L,M) \text{ is a Fredholm pair}}. \]

\sub{Regular operators.}
A (probably unbounded) operator on a Hilbert space $H$ is called regular if it is closed and densely defined.
Let $\Reg(H)$ denote the space of all regular operators on $H$. 
In this paper, we always consider $\Reg(H)$ with the \textit{graph topology}\footnote{This topology 
is strictly weaker than the so called Riesz topology,
which is induced by the embedding $\Reg(H)\hookto\B(H)$, $A\mapsto A(1+A^*A)^{-1/2}$.},
which is induced by the inclusion 
\[ \Reg(H)\hookto\Gr(H\oplus H), \quad \Op\mapsto\gr_{\Op} \] 
taking a regular operator to its graph.
The restriction of the graph topology to the subset $\B(H)$ of bounded operators 
coincides with the usual norm topology \cite[Addendum, Theorem 1]{CL}, 
so $\B(H)$ is a subspace of $\Reg(H)$.
Therefore, we have two canonical embeddings
	\begin{equation}\label{diag:BRGr}
			\B(H) \subset \Reg(H) \subset \Gr(H\oplus H).
	\end{equation}

A regular operator $\Op$ is called: 
\begin{itemize}[topsep=-2pt, itemsep=0pt, parsep=3pt, partopsep=0pt]
	\item \textit{invertible} if $\Op\colon\dom(\Op)\to H$ is bijective;
	\item \textit{Fredholm} if $\Op$ has closed range and finite-dimensional kernel and cokernel;
	\item \textit{lower semi-Fredholm} if $\Op$ has closed range and finite-dimensional kernel;
	\item \textit{upper semi-Fredholm} if $\Op$ has closed range and finite-dimensional cokernel.
\end{itemize}

\sub{Operators and their graphs.}
It is well known that properties of a regular operator $\Op$ 
can be expressed in terms of its graph $\Gamma\subset\Hh=H\oplus H$:
\begin{enumerate}[topsep=-2pt, itemsep=0pt, parsep=3pt, partopsep=0pt]
		\item The range of $\Op$ is closed if and only if the sum $\Gamma+(H\oplus 0)$ is closed in $\Hh$.
		\item $\Op$ is invertible if and only if the pair $(\Gamma,H\oplus 0)$ is transversal in $\Hh$.
		\item $\Op$ is Fredholm if and only if the pair $(\Gamma,H\oplus 0)$ is Fredholm in $\Hh$.
		\item $\Op$ is lower (resp. upper) semi-Fredholm if and only if 
		     the pair $(\Gamma,H\oplus 0)$ is lower (resp. upper) semi-Fredholm in $\Hh$.
\end{enumerate}

Indeed, identifying both $M=H\oplus 0$ and $0\oplus H$ with $H$, 
one can write the kernel and range of $\Op$ in terms of its graph:
\[ \ker(\Op) = \Gamma\cap M \txt{and} \ran \Op = (\Gamma+M)/M \subset \Hhat/M = 0\oplus H. \]
In particular, the kernel of $\Op$  has the same dimension as $\Gamma\cap M$.
The range of $\Op$ is closed if and only if $\Gamma+M$ is closed in $\Hhat$;
coincides with $H$ if and only if $\Gamma+M = \Hhat$;
and has finite codimension in $H$ if and only if $\Hh/(\Gamma+M)$ is finite-dimensional.

\section{$K^0$ index}\label{sec:K0ind}

\upskip
\sub{$K^0$ index for bounded Fredholm operators.}
Let $H$ be an infinite-dimensional Hilbert space.
By a classical result of Atiyah and J\"anich \cite{Atiyah}, 
the space $\BF(H)$ of bounded Fredholm operators (with the norm topology)
is a classifying space for the functor $K^0$.
As a set, $K^0(X)$ might be \textit{defined} as the set $[X,\BF]$ of homotopy classes of maps $X\to\BF(H)$.
The class of a continuous map $f\colon X\to\BF$ in $K^0(X)$ is called the \textit{index} of $f$:
\[ \ind(f) = [f]\in[X,\BF]=K^0(X). \]
The Kuiper's theorem implies that $[vfv\inv] = [f]$ for any continuous family of unitary operators 
$v\colon X\to\U(H)$ or, more generally, continuous family $v\colon X\to\GL(H)$ of isomorphisms.
Hence it does not matter which Hilbert space $H$ is chosen for the definition of $K^0(X)$:
one can use an arbitrary isomorphism $H'\to H$ for the identification of $[X,\BF(H')]$ with $[X,\BF(H)]$.

The group structure on $K^0(X)=[X,\BF]$ is defined as follows.
The constant map $f_x\equiv 1$ represents the neutral element of the group. 
The point-wise adjoint map $f^*\colon X\to\BF$ represents the element inverse to $[f]$.
The group operation in $[X,\BF]$ is usually defined by the composition of operators: 
$[f]+[g]=[fg]$ for $f,g\colon X\to\BF(H)$. 
Equivalently, the group operation may be defined using the direct sum of operators: 
\begin{equation*}%\label{eq:BF+}
	[f]+[g] = [f\oplus g] \in \brr{X,\BF(H\oplus H)} = K^0(X)  \quad\text{for } f,g\colon X\to\BF(H); 
\end{equation*}
this definition is more convenient for our purposes.

\begin{rem}\label{rem:K0}
The definition of $K^0(X)$ as $[X,\BF]$ works for an arbitrary space $X$, not necessarily compact.
In contrast with this, the usual topological definition of $K^0(X)$ as the Grothendieck group 
of isomorphism classes of vector bundles over $X$ is suited only for compact spaces $X$.
As was shown by Atiyah \cite{Atiyah}, for compact $X$ these two definitions are equivalent.
\end{rem}

\upskip
\sub{$K^0$ index for regular Fredholm operators.} 
Restricting the embeddings \eqref{diag:BRGr} to the subspaces of Fredholm operators, 
we obtain the canonical embeddings 
\begin{equation}\label{eq:BF-RF-GrF}
	\BF(H) \subset \RF(H) \subset \Gr_F(H\oplus H; H\oplus 0), 
\end{equation}
where $\RF(H)$ denotes the subspace of $\Reg(H)$ consisting of Fredholm operators.   

As was shown by Joachim, $\RF(H)$ is a classifying space for the functor $K^0$.
See \cite[Theorem 3.5(i)]{Jo}.
An alternative proof of this result given by the author shows that 
both embeddings in \eqref{eq:BF-RF-GrF} are homotopy equivalences.
See \cite[Theorem A]{Pr21}.
This allows to extend the definition of the index from families of bounded Fredholm operators 
to graph continuous families of regular Fredholm operators, 
and also to families of subspaces of $H\oplus H$ which are Fredholm with respect to $H\oplus 0$.

Thus, a graph continuous family $\Op=(\Op_x)$, $\xX$ of regular Fredholm operators 
parametrized by points of a topological space $X$ has an index 
\[ \ind(\Op) = [\Op] \in [X, \RF(H)] = K^0(X). \]
Similarly, a continuous family $L=(L_x)$, $\xX$ of subspaces $L_x\in\Gr_F(H\oplus H; H\oplus 0)$ 
has an index 
\[ \ind(L) = [L] \in [X, \Gr_F(H\oplus H; H\oplus 0)] = K^0(X). \]
Since \eqref{eq:BF-RF-GrF} are homotopy equivalences, 
the index of a family $\Op=(\Op_x)$ of regular Fredholm operators 
coincides with the index of the family of their graphs:
\begin{equation}\label{eq:indA}
	\ind(\Op) = \ind(\gr_{\Op}) \in K^0(X).  
\end{equation}
This observation is essential to our proof of Theorems \ref{Thm:ind0} and \ref{Thm:ind-AL}.

\sub{$K^0$ index for Fredholm pairs.}
More generally, let $f\colon X\to\Gr_F(\Hhat)$, $f(x)=(L_x,M_x)$ 
be a continuous family of Fredholm pairs in $\Hhat$.
Suppose first that all $M_x$, and thus all $L_x$, have infinite dimension and codimension
(we will discuss the general case in the next subsection).
By Proposition \ref{prop:decomp},
there is a continuous family of unitary operators $u\colon X\to\U(\Hhat)$
taking the family $M_x$ to the constant family $u_xM_x\equiv M_0$.
Identifying both $M_0$ and $M_0^{\bot}$ with an auxiliary Hilbert space $H$,
we obtain an identification of $\Hhat=M_0\oplus M_0^{\bot}$ with $H\oplus H$ taking $M_0$ to $H\oplus 0$.
The index of $f$ is then defined by the formula
\begin{equation}\label{eq:ind-LM-uL}
	\ind(f) = \ind(L,M) := \ind(uL) \in [X,\Gr_F(H\oplus H; H\oplus 0)] = K^0(X), 
\end{equation}
where $uL$ denotes the family of subspaces $u_xL_x$, that is, the map 
\[ X \to \Gr_F(\Hhat;M_0) = \Gr_F(H\oplus H; H\oplus 0), \quad x\mapsto u_xL_x. \]
This definition does not depend on the choice of $u$. 
Indeed, every two maps $u,u'\colon X\to\U(\Hhat)$ satisfying the condition $u_xM_x = H\oplus 0$
differ by a family of unitary operators $v=u'u\inv\colon X\to\U(\Hhat)$ stabilizing $H\oplus 0$, that is, 
by unitary operators of the form $v=v_0\oplus v_1$ for some $v_0,v_1\colon X\to\U(H)$.
The unitary group $\U(H)$ is contractible by Kuiper's theorem, 
so classes of $u'L = vuL$ and $uL$ in $[X,\Gr_F(\Hhat;H\oplus 0)]$ coincide.

We can now rewrite formula \eqref{eq:indA} in terms of Fredholm pairs:
\begin{equation}\label{eq:indA2}
 \ind(\Op) = \ind_{\Hh}(\gr_{\Op},H\oplus 0) \in K^0(X).
\end{equation}
(We will often use a subscript with the index: $\ind(L,M)=\ind_{\Hhat}(L,M)$,
in order to indicate the ambient space where a Fredholm pair lives.)

If $(L_x,M_x)$ is a family of Fredholm pairs consisting of subspaces of infinite dimension and codimension,
then its homotopy type is determined by its index:
two such families are homotopic as maps to $\Gr^2_F$ if and only if their indices coincide.

\sub{Stabilization.}
In the discussion above, it was supposed that the subspaces $L_x$ and $M_x$ 
have infinite dimension and codimension; in particular, $\Hhat$ itself is infinite-dimensional.
However, in our framework the space $\VV$ of boundary values may be finite-dimensional,
so we need to extend the definition of the index to a more general case.
Using stabilization, one can define the index of an arbitrary family $(L,M)$ of Fredholm pairs by the formula
\[ \ind_{\Hhat}(L,M) := \ind_{\Hhat\oplus H_0\oplus H_1}(L\oplus H_0,M\oplus H_1) \in [X,\Gr^2_F(\Hhat\oplus H_0\oplus H_1)] = K^0(X), \]
where $H_0$, $H_1$ are auxiliary infinite-dimensional Hilbert spaces.
This definition agrees with the previous one given by \eqref{eq:ind-LM-uL},
since the index \eqref{eq:ind-LM-uL} is additive with respect to the direct sums
and the index of the transversal pair $(H_0,H_1)$ in $H_0\oplus H_1$ vanishes.
Therefore, the index $\ind_{\Hhat}(L,M)\in K^0(X)$ of an arbitrary family of Fredholm pairs 
$(L_x,M_x)$ in $\Hh$ is well defined.

If the Hilbert space $\Hhat$ is \textit{finite-dimensional}, 
then the subspaces $L_x$ and $M_x$ give rise to vector bundles $L$ and $M$ over $X$. 
In this case, the index is equal to the difference of the classes of these vector bundles in $K^0(X)$: 
\[ \ind(L,M) = \ind(L,0)-\ind(M,0) = [L]-[M] \in K^0(X). \]

\upskip
\sub{Properties.} 
We will use the following properties of the index in our proofs:
\begin{enumerate}[topsep=0pt, itemsep=0pt, parsep=3pt, partopsep=0pt]
	\item Homotopy invariance: $\ind(L,M)=\ind(L',M')$ for homotopic families of Fredholm pairs 
	      $(L_x,M_x)$ and $(L'_x,M'_x)$ in $\Hh$.
	\item Additivity with respect to orthogonal sums: 
	      \[ \ind_{\Hh\oplus\Hh'}(L\oplus L',M\oplus M') = \ind_{\Hh}(L,M)+\ind_{\Hh'}(L',M'), \]
				for families $(L_x,M_x)$ of Fredholm pairs in $\Hh$ and $(L'_x,M'_x)$ in $\Hh'$.
	\item Invariance under isomorphisms: $\ind(L,M)=\ind(uL,uM)$ for a continuous family 
	      of unitary operators $u\colon X\to\U(\Hh)$.
	\item Vanishing on families of transversal pairs.
\end{enumerate}

\medskip
The first two properties hold by the definition.
The third property follows from the fact that \eqref{eq:ind-LM-uL} does not depend on the choice of $u$, 
see discussion below \eqref{eq:ind-LM-uL}.
Let us explain the last property.
A pair $(L,M)$ is transversal if and only if $L$ is the graph of a bounded operator $T_L\colon M^{\bot}\to M$.
Using stabilization if necessary,
we can suppose that the family $M$ is already transformed into the constant family $M_x\equiv H\oplus 0$.
A linear deformation of $T_L$ to zero operator leads to the deformation of the family $(L,M)$
to the constant family $(M^{\bot},M) = (0\oplus H, H\oplus 0)$, whose index vanishes.
Therefore, the index of $(L,M)$ also vanishes.

\sub{Hilbert space depending on parameter.}
More generally, the Hilbert space $\Hh=\Hh_x$ may also depend on $x$.
Let $\H$ be a locally trivial Hilbert bundle over a paracompact space $X$ 
with the fiber $\Hh$ and the structure group $\U(\Hh)$.
Then the index of a family $(L_x,M_x)\in\Gr_F(\H_x)$ of Fredholm pairs living in the fibers of $\H$
is well defined.
Indeed, if $\dim H$ is infinite, then the bundle $\H$ is trivial by Kuiper's theorem. 
A trivialization $\H\cong \Hh\times X$ takes $(L,M)$ to a family of Fredholm pairs in a fixed Hilbert space $\Hh$,
and property (2) above implies that the index of thus obtained family 
is independent of the choice of a trivialization.
The case of finite-dimensional $\Hh$ is obtained from the infinite-dimensional case using stabilization.

\section{Extension pairs of operators}\label{sec:ext1}

\upskip
\sub{Cauchy data space.}
The \textit{(abstract) Cauchy data space} of an extension pair $\A = (\Ami,\Ama)$
is defined as the trace of the kernel of $\Ama$,
\[ \CDS = \CDS_{\A} = \gamma(\ker\Ama)\subset\VV_{\A}. \]

\begin{prop}\label{prop:C-closed}
If the range of $\Ami$ is closed, then $\CDS_{\A}$ is a closed subspace of $\VV_{\A}$.
\end{prop}

\proof
Consider the diagram 
\begin{equation}\label{diag:AKC}
	\begin{tikzcd}
	  \ran\Ami \hookr{d} & \dom\Ami+\ker\Ama \arrow[swap]{l}{\Ama} \arrow{r}{\gamma} \hookr{d} & \CDS \hookr{d}
		\\
		\ran\Ama & \dom\Ama \arrow[swap]{l}{\Ama} \arrow{r}{\gamma} & \VV
	\end{tikzcd}
\end{equation}
(here all domains are equipped with the graph norm, as usual).
Since $\gamma$ is a quotient map, $\CDS$ is closed in $\VV$ if and only if its inverse image 
$\gamma\inv(\CDS) = \dom\Ami+\ker\Ama$ is closed in $\dom\Ama$. 
But $\dom\Ami+\ker\Ama$ is the inverse image of $\ran\Ami$ under the bounded operator 
$\Ama\colon\dom\Ama\to H$, so it is closed if $\ran\Ami$ is closed.
\endproof

The following lemma provides somewhat converse statement to the last sentence of the proof above.

\begin{lem}\label{lem:ran2}
	Let $\Ami\subset\Ama$ be an extension pair such that the range of $\Ama$ is closed.
	If $\dom\Ami+\ker\Ama$ is closed in $\dom\Ama$, then the range of $\Ami$ is also closed.
\end{lem}

\proof
Consider the left commutative square on Diagram \eqref{diag:AKC}.
The operator 
\begin{equation}\label{eq:Ama}
	\Ama\colon\dom\Ama\to\ran\Ama 
\end{equation}
is bounded with respect to the graph norm on $\dom\Ama$.
It induces the bijection 
\[ \dom\Ama/\ker\Ama\to\ran\Ama \] 
between Hilbert spaces.
By the Closed Graph Theorem, this bijection is a topological isomorphism, so \eqref{eq:Ama} is a quotient map.
The inverse image of $\ran\Ami$ under this map is $\dom\Ami+\ker\Ama$, 
which is closed in $\dom\Ama$ by the assumption of the lemma.
It follows that $\ran\Ami$ is closed in $\ran\Ama$ and thus in $H$.
\endproof

\sub{Properties of realizations.} 
Let $\A = (\Ami,\Ama)$ be an arbitrary extension pair and $L$ be a subspace of $\VV$.
Since $\ker\A_L = \dom\A_L\cap\ker\Ama$ and $\dom\Ami\subset\dom\A_L$, we get 
\[ \gamma(\ker\A_L) = \gamma(\dom\A_L)\cap\gamma(\ker\Ama) = L\cap\CDS. \]
Since the kernel of the restriction $\gamma\colon\ker\A_L\to L\cap\CDS$ is equal to 
$\ker\A_L\cap\dom\Ami = \ker\Ami$, we obtain a canonical exact sequence of vector spaces
\begin{equation}\label{eq:ex-sec1}
		0 \to \ker\Ami \to \ker\A_L \stackrel{\gamma}{\longrightarrow} L\cap\CDS \to 0. 				
\end{equation}

Suppose that $\ran\A_L$ is closed.
Replacing the top row in Diagram \eqref{diag:AKC} by the row
\[ \ran\A_L \stackrel{\Ama}{\longleftarrow} \dom\A_L+\ker\Ama \stackrel{\gamma}{\longrightarrow} L+\CDS \]
and reasoning as in the proof of Proposition \ref{prop:C-closed}, we see that $L+\CDS$ is closed. 

Conversely, suppose that $L+\CDS$ is closed. 
Then the subspace $\dom\A_L+\ker\Ama = \gamma\inv(L+\CDS)$ is closed in $\dom\Ama$.
If additionally $\ran\Ama$ is closed, then Lemma \ref{lem:ran2} applied to the pair $\A_L\subset\Ama$
implies that $\ran\A_L$ is closed.
Moreover, we have a canonical exact sequence of Hilbert spaces
\begin{equation}\label{eq:ex-sec2}
		0 \to \coker\Ama \to \coker\A_L \to \VV/(L+\CDS) \to 0. 
\end{equation}
Indeed, for the boundary condition $N=L+\CDS$, the domain of the corresponding realization is  
$\dom\A_N=\dom\A_L+\ker\Ama$, so $\ker\A_N=\ker\Ama$ and $\ran\A_N=\ran\A_L$.
Hence
\[ \ran\Ama/\ran\A_L = \ran\Ama/\ran\A_N \cong \dom\Ama/\dom\A_N \cong \VV/N = \VV/(L+\CDS). \]
This provides the canonical isomorphism $\ran\Ama/\ran\A_L \cong \VV/(L+\CDS)$,
from which \eqref{eq:ex-sec2} follows.

\sub{Fredholm extension pairs.}
We say that an extension pair $\A = (\Ami,\Ama)$ is \textit{Fredholm} 
if the minimal operator $\Ami$ is lower semi-Fredholm and the maximal operator $\Ama$ is upper semi-Fredholm.
Similarly, we say that $\A$ is \textit{invertible} if 
$\Ami$ is injective and has closed range and $\Ama$ is surjective.
For $\Ami=\Ama$ these notions coincide with the usual Fredholmness and invertibility.

Applying conclusions of  the previous subsection to a Fredholm extension pair $\A$,
we see that the range of $\A_L$ is closed if and only if $L+\CDS$ is closed.
Moreover, 
\[ \dim\ker\A_L = \dim\ker\Ami+\dim(L\cap\CDS) \text{ and } \dim\coker\A_L = \dim\coker\Ama+\codim(L+\CDS). \]
This provides the following criterion.

\begin{prop}\label{prop:AL-Fr}
  Let $\A = (\Ami,\Ama)$ be a Fredholm (resp. invertible) extension pair 
	and $L$ be a closed subspace of $\VV$.
	Then the realization $\A_L$ is Fredholm (resp. invertible) if and only if 
	$(L,\CDS)$ is a Fredholm (resp. transversal) pair of subspaces of $\VV$.
\end{prop}

The following property provides an additional justification for our terminology concerning extension pairs;
it is not used in this paper.

\begin{prop}\label{prop:Amima-tF}
Let $\A = (\Ami,\Ama)$ be an extension pair.
Then $\A$ is Fredholm (resp. invertible) if and only if $\A$ admits Fredholm (resp. invertible) realizations.
\end{prop}

\proof
($1\Rightarrow 2$) 
Let $\A$ be Fredholm (resp. invertible).
Choosing $L\subset\VV$ such that the pair $(L,\CDS)$ is transversal (for example, $L=\VV\ominus\CDS$),
we obtain a Fredholm (resp. invertible) realization $\A_L$.

($2\Rightarrow 1$) 
Let $\At$ be an invertible realization of $\A$.
Then $\ran\Ama\supset\ran\At=H$ and $\ker\Ami\subset\ker\At=0$.
In addition, $\At$ induces a topological isomorphism $\dom\At\to H$ 
(with respect to the graph norm on $\dom\At$)
and thus takes a closed subspace $\dom\Ami$ of $\dom\At$ to a closed subspace $\ran\Ami$ of $H$.
Therefore, the extension pair $\A$ is invertible.

Let now $\At$ be a Fredholm realization of $\A$. 
Then the range of $\At$ has finite algebraic codimension, 
so $\ran\Ama\supset\ran\At$ also has finite algebraic codimension.
In other words, $\Ama$ is upper semi-Fredholm.
The kernel of $\Ami$ is contained in $\ker\At$ and thus has finite dimension.
Since $\ker\At$ is finite-dimensional and $\dom\Ami$ is closed in $\dom\At$,
the sum $\dom\Ami+\ker\At$ is also closed in $\dom\At$.
Applying Lemma \ref{lem:ran2} to the pair $\Ami\subset\At$, we see that $\ran\Ami$ is closed.
Therefore, $\Ami$ is lower semi-Fredholm and the extension pair $\A$ is Fredholm.
\endproof

\section{Pairs $\Gamma\subset\Gamma'$}\label{sec:Gamma}

As we saw in Section \ref{sec:GrF},
properties of operators in $H$ can be expressed in terms of their graphs in $\Hh=H\oplus H$.
In this section we generalize constructions and results of the previous section to arbitrary subspaces 
of a Hilbert space $\Hh$.
Some proofs from the previous section become even simpler in this more general framework; 
we explain in the end of this section how they are related.

\sub{Intermediate subspaces.}
Let $\Gamma\subset\Gamma'$ be closed subspaces of a Hilbert space $\Hh$,
$\VV=\Gamma'/\Gamma$ be the quotient Hilbert space,
and $\gamma\colon \Gamma'\to\VV$ be the quotient map.

The inverse image $\gamma\inv(L)$ of a closed subspace $L\subset\VV$ 
is a closed subspace in $\Gamma'$ containing $\Gamma$.
Considering $\gamma\inv(L)$ as a subspace of $\Hhat$, we obtain the map
\begin{equation}\label{eq:gammi}
	\gammi \colon\Gr(\VV)\hookto\Gr(\Hh), \quad L\mapsto \gamma\inv(L). 
\end{equation}
We will often omit parentheses in $\gamma\inv(L)$ in order to simplify notations.

Conversely, every closed subspace $\Lt$ of $\Hhat$ such that $\Gamma\subset\Lt\subset \Gamma'$ 
has the form $\gammi L$ for $L=\gamma(\Lt)\subset\VV$.
Therefore, closed subspaces $\Lt$ lying between $\Gamma$ and $\Gamma'$ 
are parameterized by the Grassmanian $\Gr(\VV)$.

\begin{prop}\label{prop:gamma*}
	Inclusion \eqref{eq:gammi} is an embedding (that is, a homeomorphism on the image).
\end{prop}

\proof
Identifying $\VV$ with the subspace $\Gamma'\ominus \Gamma$ of $\Hh$, 
we can write $\gammi L = L\oplus\Gamma$. 
The orthogonal projection onto $\gammi L$ is equal to $P_L+P_{\Gamma}$ and thus depends continuously on $L$. 
In other words, the inclusion \eqref{eq:gammi} is continuous.
Conversely, $P_L = P_{\gammi L}-P_{\Gamma}$ depends continuously on $\gammi L$, so \eqref{eq:gammi} is an embedding. 
\endproof

\sub{Push-forward.} 
In addition to \eqref{eq:gammi}, there is another map, which we denote by $\gamma_!$, 
going in the opposite direction and taking a subspace $M$ of $\Hh$ to a subspace of $\VV$:
\begin{equation}\label{eq:gamma!}
  \gamma_!M := \gamma(M\cap \Gamma') \subset \VV.
\end{equation}
Clearly, $\gamma_!\gammi =\Id$.
If $M\subset\Gamma'$, then $\gamma_!M=\gamma(M)$, the image of $M$ under the map $\gamma\colon \Gamma'\to\VV$.
The subspace $\gamma_!M$ can be considered as the \q{trace}, or approximation of $M$ in $\VV$.

\begin{prop}\label{prop:!closed}
	If $\Gamma+M$ is closed in $\Hh$, then $\gamma_!M$ is closed in $\VV$.
\end{prop}

\proof
If $\Gamma+M$ is closed in $\Hh$, then 
\[ \gamma\inv(\gamma_!M) = (\Gamma+M)\cap \Gamma' \]
is closed in $\Gamma'$. 
Since $\gamma$ is a quotient map, $\gamma_!M$ is closed in $\Gamma'/\Gamma=\VV$.
\endproof

\sub{Relation between pairs $(\gammi L,M)$ and $(L,\gamma_!M)$.}
Let $M$ be a subspace of $\Hh$ and $L$ be a subspace of $\VV$.
Denoting $M'=M\cap\Gamma'$, we get 
\[ \gamma_!\br{\gammi L\cap M} = \gamma\br{\gamma\inv L\cap M'} = L\cap\gamma(M') = L\cap\gamma_!M. \]
The kernel of the surjective map $\gamma\colon \gamma\inv(L)\cap M' \to L\cap\gamma(M')$ is
\[ \Gamma\cap\br{\gamma\inv(L)\cap M'} = \Gamma\cap M' = \Gamma\cap M, \] 
so we obtain a canonical exact sequence of vector spaces
\begin{equation}\label{eq:exseq-LM1}
	0 \to \Gamma\cap M \to \gammi L\cap M \stackrel{\gamma}{\longrightarrow} L\cap\gamma_!M \to 0. 
\end{equation}
Similarly, 
\[ \gamma_!(\gammi L+M) = \gamma(\gammi L+M') = L+\gamma(M') = L+\gamma_!M. \]
Since $\gamma$ is a quotient map and $\gammi L+M'$ contains $\ker\gamma=\Gamma$, we see that 
\begin{equation*}%\label{eq:gamma-LM3}
(\gammi L+M)\cap\Gamma' = \gammi L+M' \text{ is closed in } \Gamma' 
   \txt{ if and only if } L+\gamma_!M \text{ is closed in } \VV.
\end{equation*}
If $\gammi L+M$ is closed in $\Hh$, then $(\gammi L+M)\cap\Gamma'$ is closed in $\Gamma'$ and thus 
$L+\gamma_!M$ is closed in $\VV$.

Conversely, suppose that $L+\gamma_!M$ is closed in $\VV$. 
Then $\gammi L+M'$ is closed in $\Gamma'$
and thus the direct sum $(\gammi L+M')\dotplus(M\ominus M') = \gammi L+M$ 
is closed in the direct sum $\Gamma'\dotplus(M\ominus M') = \Gamma'+M$.
Suppose additionally that $\Gamma'+M$ is closed in $\Hh$. 
Then $\gammi L+M$ is also closed in $\Hh$ and we get
\[ (\Gamma'+M)/(\gammi L+M) \;\cong\; \Gamma'/(\gammi L+M') \;=\; 
   \Gamma'/\gamma\inv(L+\gamma_!M) \cong \VV/(L+\gamma_!M), \]
where both isomorphisms are canonical.
Therefore, we have a canonical exact sequence of Hilbert spaces
\begin{equation}\label{eq:exseq-LM2}
 0 \to \Hh/(\Gamma'+M) \to \Hh/(\gammi L+M) \to \VV/(L+\gamma_!M) \to 0. 
\end{equation}

\sub{Transversal and Fredholm pairs.}
We will say that a closed subspace $M$ of $\Hh$ is \textit{Fredholm} with respect to the pair $(\Gamma,\Gamma')$
if $(M,\Gamma)$ is lower semi-Fredholm and $(M,\Gamma')$ is upper semi-Fredholm.
We will say that $M$ is \textit{transversal} to $(\Gamma,\Gamma')$ if
\begin{equation*}%\label{eq:MN-trans}
	M+\Gamma \text{ is closed}, \quad M\cap \Gamma = 0, \quad M+\Gamma'=\Hh.  
\end{equation*}
For $\Gamma=\Gamma'$ these notions coincide with the usual Fredholmness and transversality 
of the pair $(M,\Gamma)$.

Applying conclusions of the previous subsection to the particular case of Fredholm or transversal $M$,
we obtain the following result.

\begin{prop}\label{cor:gamma-LM-tF}
Let $M$ be Fredholm (resp. transversal) with respect to $(\Gamma,\Gamma')$, 
and let $L$ be a closed subspace of $\VV$. 
Then the pair $(\gammi L,M)$ is Fredholm (resp. transversal) 
if and only if the pair $(L,\gamma_!M)$ is Fredholm (resp. transversal). 
\end{prop}

The following property provides an additional justification for our terminology;
it is not used in this paper. 

\begin{prop}\label{prop:MGG'L}
	Let $M$ and $\Gamma\subset\Gamma'$ be closed subspaces of $\Hh$.
	Then the following two conditions are equivalent:
	\begin{enumerate}[topsep=-2pt, itemsep=0pt, parsep=3pt, partopsep=0pt]
		\item $M$ is Fredholm with respect to (resp. transversal to) $(\Gamma,\Gamma')$.
		\item There exists a closed subspace of $\Hh$ lying between $\Gamma$ and $\Gamma'$ 
		and Fredholm with respect to (resp. transversal to) $M$.
	\end{enumerate}
\end{prop}

\proof
($1\Rightarrow 2$) 
Taking $L=\VV\ominus\gamma_!M$, we obtain a desired subspace $\gammi L$ of $\Hh$.

($2\Rightarrow 1$) 
Let $\Gamma\subset N\subset\Gamma'$.
Suppose first that the pair $(M,N)$ is transversal.
Then $M\cap\Gamma = 0$ and $M+\Gamma'=\Hh$.
The addition map $M\oplus N\to\Hh$, $\xi\oplus\eta\mapsto\xi+\eta$ is an isomorphism
which takes $M\oplus\Gamma$ to $M+\Gamma$.
Since $M\oplus\Gamma$ is closed in $M\oplus N$, the sum $M+\Gamma$ is closed in $\Hh$.
Hence $M$ is transversal to $(\Gamma,\Gamma')$.

Suppose now that $(M,N)$ is Fredholm.
Then $M+N$ admits a complement subspace of finite dimension in $\Hhat$,
so $M+\Gamma'\supset M+N$ also admits a complement subspace of finite dimension,
hence $(M,\Gamma')$ is upper semi-Fredholm.
In particular, $M+\Gamma'$ is closed.
Replacing $M$ by $M'=M\ominus(M\cap N)$, we obtain the pair $(M',N)$ 
which is transversal in the Hilbert space $M+\Gamma'$.
As was shown above, this implies that $M'+\Gamma$ is closed in $M+\Gamma'$ and thus in $\Hh$.
Since $M+\Gamma=(M'+\Gamma)+M\cap N$ and $M\cap N$ is finite-dimensional, 
the subspace $M+\Gamma$ is also closed in $\Hh$. 
Since $M\cap\Gamma\subset M\cap N$ is finite-dimensional, the pair $(M,\Gamma)$ is lower semi-Fredholm.
\endproof

\sub{Extension pairs, realizations, and graphs.}
Let us explain how the previous section fits into the framework of this section.

Let $\A = (\Ami,\Ama)$ be an extension pair of operators in $H$.
It leads to the pair $\Gamma\subset\Gamma'$ of subspaces of $\Hhat=H\oplus H$,
with $\Gamma$ being the graph of $\Ami$ and $\Gamma'$ being the graph of $\Ama$.
The space $\VV_{\A}=\dom\Ama/\dom\Ami$ of abstract boundary values 
is identified with the quotient space $\Gamma'/\Gamma$, 
and the trace map 
\[ \gamma\colon\dom\Ama\to\dom\Ama/\dom\Ami=\VV_{\A} \]
is identified with the quotient map $\Gamma'\to\Gamma'/\Gamma$.

Expressing properties of operators in terms of their graphs as explained in Section \ref{sec:GrF}, we see that 
$\A$ is Fredholm if and only if $H \oplus 0$ is Fredholm with respect to $(\Gamma,\Gamma')$ in $\Hh$;
$\A$ is invertible if and only if $H \oplus 0$ is transversal to $(\Gamma,\Gamma')$ in $\Hh$.

Let $\At$ be a realization of $\A$.
Then the graph $\gr_{\At}$ of $\At$ is a closed subspace of $\Hhat$ 
such that $\Gamma\subset\gr_{\At}\subset\Gamma'$.
Such subspaces are in one-to-one correspondence with closed subspaces of $\VV=\Gamma'/\Gamma$.
A boundary condition $L\subset\VV$ gives rise to the realization $\A_L$ with the graph $\gammi L$.
Conversely, a realization $\At$ is given by the boundary condition 
$L = \gamma(\gr_{\At}) = \gamma_!\gr_{\At}$.

The Cauchy data space $\CDS_{\A}$ can be written in terms of the push-forward map $\gamma_!$.
Indeed, the identification of $\dom(\Ama)$ with $\Gamma'$ takes $\ker\Ama$ 
to the intersection $\Gamma'\cap(H\oplus 0)$ 
and the image $\gamma(\ker\Ama)$ to $\gamma(\Gamma'\cap(H\oplus 0)) = \gamma_!(H\oplus 0)$.
Therefore, 
\begin{equation*}%\label{eq:C-gamma!}
	\CDS = \gamma_!(H\oplus 0). 
\end{equation*}
In such a way, all the results of this section can be carried over to extension pairs and realizations.
Applying Propositions \ref{prop:!closed}, \ref{cor:gamma-LM-tF}, and \ref{prop:MGG'L} 
to $\Gamma$, $\Gamma'$, $M=H\oplus 0$, and $\CDS = \gamma_!M$,
one obtains Propositions \ref{prop:C-closed}, \ref{prop:AL-Fr}, and \ref{prop:Amima-tF}. 
Similarly, one obtains exact sequences \eqref{eq:ex-sec1} and \eqref{eq:ex-sec2} from
exact sequences \eqref{eq:exseq-LM1} and \eqref{eq:exseq-LM2}.

\section{Push-forward for transversal subspaces}\label{sec:TN}

If $M$ is Fredholm with respect to $(\Gamma,\Gamma')$, 
then $\gamma_! M$ is a closed subspace of $\VV$ by Proposition \ref{prop:!closed}.
Thus $\gamma_!$ determines a map from the Fredholm, with respect to $(\Gamma,\Gamma')$, 
part of $\Gr(\Hh)$ to $\Gr(\VV)$.
In contrast with $\gammi$, this map is \textit{not} continuous.
However, it behaves better on a smaller part of $\Gr(\Hh)$,
\[ \T = \sett{M\in\Gr(\Hh)}{M \text{ is transversal to } (\Gamma,\Gamma')}. \]

\begin{prop}\label{prop:tr-tr}
The restriction $\gamma_!\colon\T\to\Gr(\VV)$ of $\gamma_!$ to $\T$ is continuous.
\end{prop}

\proof
By Proposition \ref{prop:cap-cup-cont}, the map $\T\to\Gr(\Gamma')$, $M\mapsto M\cap\Gamma'$ is continuous.
The image of this map consists of subspaces $M'\subset\Gamma'$ such that $M'\cap\Gamma=0$.
Applying Proposition \ref{prop:cap-cup-cont} one more time, 
we see that the map $M'\mapsto M'+\Gamma$ is continuous,
and thus also $M'\mapsto (M'+\Gamma)/\Gamma = \gamma_!M$ is continuous.
\endproof

\medskip
For a closed subspace $N$ of $\VV$, let $\T_N$ be the fiber of $\T$ over $N$: 
\[ \T_N = \sett{M\in\T}{\gamma_!M = N}. \] 
Identifying $\VV$ with the subspace $\Gamma'\ominus \Gamma$ of $\Hh$, 
we obtain a distinguished element 
\[ \bar{N}=N\oplus(\Hh\ominus\Gamma') \in\T_N. \] 
The map $N\mapsto \bar{N}$ is a continuous section of the bundle $\gamma_!\colon\T\to\Gr(\VV)$.

\begin{prop}\label{prop:TN}
	The space $\T_N$ has a canonical structure of a Banach space,
	with $\bar{N}$ being its zero element.
\end{prop}

\proof
Denote $E=\VV\ominus N$ and $F=\Hh\ominus\Gamma'$,
and consider the orthogonal decomposition
\begin{equation*}%\label{eq:Hh-dec}
	\Hh = (E\oplus\Gamma)\oplus(N\oplus F) = (E\oplus\Gamma)\oplus\bar{N}.
\end{equation*}
Suppose that $M\in\T_N$.
Since $E$ is transversal to $\gamma_!M = N$ in $\VV$ and $M$ is transversal to $(\Gamma,\Gamma')$ in $\Hh$, 
it follows from Proposition \ref{cor:gamma-LM-tF} that $M$ is transversal to $\gammi E=E\oplus\Gamma$.
Hence $M$ is the graph of a bounded map 
\begin{equation}\label{eq:mu}
	\mu\colon N\oplus F \to E\oplus\Gamma. 
\end{equation}
The intersection $M\cap\Gamma'$ is the graph of the restriction $\restr{\mu}{N}\colon N \to E\oplus\Gamma$.
The last map can be written as the direct sum 
\[ \restr{\mu}{N} = \restr{P_E\mu}{N}\oplus\restr{P_{\Gamma}\mu}{N}, \]
where $P_E$ and $P_{\Gamma}$ denote the orthogonal projections $E\oplus\Gamma\to E$ and $E\oplus\Gamma\to\Gamma$.
In these notations, $\gamma_!M\subset\VV = N\oplus E$ 
is the graph of the composition $\restr{P_E\mu}{N}\colon N\to E$.
Therefore, $\gamma_!$ takes the graph of \eqref{eq:mu} to $N$ if and only if 
\begin{equation}\label{eq:phi-NE0}
 \restr{P_E\mu}{N} = 0.
\end{equation}
Conversely, suppose that $M$ is the graph of a bounded operator \eqref{eq:mu}.
Then $M$ is transversal to $E\oplus\Gamma$ and thus also transversal to $(\Gamma,\Gamma')$.
If additionally $\mu$ satisfies \eqref{eq:phi-NE0}, 
then $\gamma_!M = N$ and thus $M\in\T_N$.

Therefore, the construction above provides the canonical homeomorphism between $\T_N$ and a Banach space, 
namely the closed linear subspace of 
$\B(N\oplus F, E\oplus\Gamma)$ consisting of operators $\mu$ satisfying linear condition \eqref{eq:phi-NE0}.
This homeomorphism takes the operator $\mu=0$ to $\bar{N}=N\oplus F$. 
\endproof

\section{Index formula for Fredholm pairs}\label{sec:ind-Gr}

\upskip
\sub{Proof of Theorem \ref{Thm:ind-Gr1}.}
Let $\CDS=\gamma_!M$ and $F=\Hh\ominus\Gamma'$.
By Proposition \ref{prop:TN}, $\T_{\CDS}$ is path connected,
so there is a path $\M\colon[0,1]\to\T_{\CDS}$ connecting $\M_0=M$ and $\M_1 = \CDS\oplus F$.
By the definition of $\T_{\CDS}$, the subspace $\M_t$ is transversal to the pair $(\Gamma,\Gamma')$ 
for every $t\in[0,1]$.
Since all the pairs $(L_x,\gamma_!\M_t) = (L_x,\CDS)$ are Fredholm, 
the pairs $(\gammi L_x,\M_t)$ are also Fredholm.
By the homotopy invariance of the index, 
\[ \ind(\gammi L,\, \M_0) = \ind(\gammi L,\, \M_1), \txt{that is,} 
    \ind(\gammi L,\, M) = \ind(\gammi L,\, \CDS\oplus F). \]
The Fredholm pair $(\gammi L_x,\, \CDS\oplus F) = (L_x\oplus \Gamma, \, \CDS\oplus F)$
is the orthogonal sum of the Fredholm pairs $(L_x,\CDS)$ in $\VV$
and $(\Gamma,F)$ in $\Hh\ominus\VV = \Gamma\oplus F$.
The second pair $(\Gamma,\, F)$ is transversal, so its index vanishes.
The additivity of the index implies 
\[ \ind_{\Hhat}(\gammi L,\, \CDS\oplus F) = \ind_{\VV}(L,\,\CDS) + \ind_{\Gamma\oplus F}(\Gamma,\,F) 
   = \ind_{\VV}(L,\,\CDS). \]
Taking this together, we obtain
	\begin{equation*}%\label{eq:thm:ind-Gr1}
		\ind_{\Hhat}(\gammi L,M) = \ind_{\VV}(L,\gamma_!M) \in K^0(X). \quad\square
	\end{equation*}

\sub{Adding parameters.}
More generally, not only $L$ but also $M$, $\Gamma$, and $\Gamma'$ may depend on a parameter $x\in X$.
Then $\VV_x = \Gamma'_x/\Gamma_x$ and $\gamma_x\colon\Gamma'_x\to\VV_x$ also depend on $x$.
We will say that a family has a property (say, is transversal or Fredholm) if each its member has this property.

Let $\Gamma_x\subset\Gamma'_x$ be continuous families of closed subspaces of $\Hh$ 
parametrized by points of a topological space $X$.
Then $\VV_x$ identified with the subspace $\Gamma'_x\ominus\Gamma_x$ of $\Hhat$ depends continuously on $x$.

If $L_x$ is a continuous family of closed subspaces of $\VV_x$, 
then the family $\gammi L$  is continuous by the same reasoning as in the proof of Proposition \ref{prop:gamma*}.

If $M=(M_x)$ is a continuous family of closed subspaces of $\Hh$ transversal to $(\Gamma,\Gamma')$, 
then the family $\gamma_! M$ is continuous by the same reasoning as in the proof of Proposition \ref{prop:tr-tr}.

\begin{thm}\label{thm:ind-Gr}
Let $M_x\subset \Hh$, $\Gamma_x\subset\Gamma'_x\subset \Hh$, and $L_x\subset\VV_x$
be continuous families of closed subspaces.
	Suppose that $M$ is transversal to $(\Gamma,\Gamma')$	and Fredholm with respect to $\gammi L$.
	Then 
	\begin{equation}\label{eq:ind-Gr}
		\ind_{\Hhat}(\gammi L,M) = \ind_{\VV}(L,\gamma_! M) \in K^0(X). 
	\end{equation}
\end{thm}

\proof
Let $\M\colon[0,1]\times X \to\T_N$ be the linear homotopy between the families $\M_0=M$ 
and $\M_1 = \gamma_! M\oplus (\Hh\ominus\Gamma')$ provided by Proposition \ref{prop:TN}.
Since $\M_{0,x}$ and $\M_{1,x}$ depend continuously on $x$, 
the subspace $\M_{t,x}$ depends continuously on $(t,x)$.
The rest of the proof is completely similar to the proof of Theorem \ref{Thm:ind-Gr1}.
\endproof

\begin{rem}\label{rem:ind-Hb}
	One can also consider a family of Hilbert spaces $\Hh_x$ depending on a parameter $x$,
	namely, the fibers of a locally trivial Hilbert bundle over $X$. 
	Since our construction of a Banach space structure in Proposition \ref{prop:TN} is canonical, 
	it works in the fibers of a Hilbert bundle as well,
	and the same reasoning as above provides a homotopy $\M$ and the equality 
	$\ind(\gammi L,M) = \ind(L,\gamma_! M)$.
\end{rem}

\upskip
\sub{Replacing transversality by constant dimension condition.}
One can generalize Theorem \ref{thm:ind-Gr} further,
replacing the assumptions $M_x\cap\Gamma_x = 0$ and $M_x+\Gamma'_x=\Hh$  
by the assumption of locally constant dimensions of the finite-dimensional spaces 
\begin{equation}\label{eq:KK'}
	K_x=M_x\cap\Gamma_x \txt{and} K'_x=\Hh\ominus(M_x+\Gamma'_x). 
\end{equation}
Under this assumption, the subspaces $K_x$ and $K'_x$ depend continuously on $x$ 
and thus give rise to locally trivial vector bundles $K$ and $K'$ over $X$.
In particular, the classes $[K],[K']\in K^0(X)$ are well defined.
The following theorem is not used in the rest of the paper; we give it here for completeness only.

\begin{thm}\label{thm:ind-Gr-KK'}
	Suppose that $M_x$ is Fredholm with respect to $\gammi L_x$ and the dimensions of the subspaces 
	\eqref{eq:KK'} are locally constant.
	Then  
	\begin{equation}\label{eq:ind-Gr-KK'}
		\ind_{\Hh}(\gammi L,M) = \ind_{\VV}(L,\gamma_! M) + [K] - [K'] \in K^0(X), 
	\end{equation}
	where $K$ and $K'$ are vector bundles over $X$ defined by \eqref{eq:KK'}.
\end{thm}

\proof
The Fredholm pair $(\gammi L,M)$ in $\Hh$ is the direct sum of 
the Fredholm pair $(\gammi L\ominus K,M\ominus K)$ in $\Hh\ominus(K\oplus K')$
and the Fredholm pair $(K,K)$ in $K\oplus K'$, so additivity of the index implies.
\begin{equation}\label{eq:ind-Gr-KK'+}
	\ind_{\Hh}(\gammi L,M) = \ind_{\Hh\ominus(K\oplus K')}(\gammi L\ominus K,M\ominus K)
	  + \ind_{K\oplus K'}(K,K). 
\end{equation}
Applying Theorem \ref{thm:ind-Gr} to $M\ominus K$, $\Gamma\ominus K$, $\Gamma'\ominus K$, 
and $\Hh\ominus(K\oplus K')$ in place of $M$, $\Gamma$, $\Gamma'$, and $\Hh$ 
and using the identity $(\Gamma'\ominus K)\ominus(\Gamma\ominus K) = \Gamma'\ominus\Gamma=\VV$,
we find the first summand in \eqref{eq:ind-Gr-KK'+}:
\[ \ind_{\Hh\ominus(K\oplus K')}(\gammi L\ominus K,M\ominus K) = \ind_{\VV}(L,\gamma_! M), \]
The second summand in \eqref{eq:ind-Gr-KK'+} is given by
\[ \ind_{K\oplus K'}(K,K) = [K\cap K] - [(K\oplus K')\ominus(K+K)] = [K] - [K']. \]
Summing up these two indices, we obtain \eqref{eq:ind-Gr-KK'}.
\endproof

\section{Index formula for Fredholm realizations}\label{sec:ind-ext}

\upskip
\sub{Continuity.}
Recall that $\Reg(H)$ denotes the space of regular operators on $H$ equipped with the graph topology.

\begin{prop}\label{prop:Gr-Reg-cont}
The natural inclusion $\Gr(\VV)\hookto\Reg(H)$, $L\mapsto\A_L$ is an embedding. 
\end{prop}

\proof
Consider the sequence of two inclusions
\[ \Gr(\VV)\hookto\Reg(H)\hookto\Gr(H\oplus H), \] 
where the first map is given by $L\mapsto	A_L$ and the second one by $\At\mapsto\gr_{\At}$.
Their composition $\Gr(\VV)\hookto\Gr(H\oplus H)$ is given by the formula 
$L\mapsto \gammi L$ and is an embedding by Proposition \ref{prop:gamma*}.
The second inclusion is an embedding by the definition of the graph topology.
Therefore, the first inclusion is also an embedding.
\endproof

\sub{Fredholm boundary conditions.}
Let $\A = (\Ami,\Ama)$ be a Fredholm extension pair.
The boundary condition $L\subset\VV_{\A}$ is called \textit{Fredholm} 
if the realization $\A_L$ is Fredholm (equivalently, the pair $(L,\CDS)$ is Fredholm).
By Proposition \ref{prop:AL-Fr}, we have the natural inclusion
 \begin{equation*}%\label{eq:GrF-RF}
	 \GrF(\VV; \CDS)\hookto\RF(H), \quad L\mapsto \A_L,
 \end{equation*}
which is an embedding by Proposition \ref{prop:Gr-Reg-cont}.
In particular, a continuous family $L=(L_x)$ of Fredholm boundary conditions 
leads to a graph continuous family of Fredholm realizations $\A_{L_x}$, 
which we will denote by $\A_{L}$.

In the following proofs, we will use notations introduced in the previous sections.
In particular, $\Gamma$ and $\Gamma'$ denote the graphs of $\Ami$ and $\Ama$.

\sub{Proof of Theorem \ref{Thm:ind0}.}
The graph of $\A_{L_x}$ coincides with $\gammi L_x$, so \eqref{eq:indA2} takes the form
\[ \ind\A_L = \ind_{\Hh}(\gammi L,H\oplus 0) \in K^0(X). \]
Since $\gamma_!(H\oplus 0) = \CDS$, Theorem \ref{Thm:ind-Gr1} implies 
\[ \ind_{\Hh}(\gammi L,H\oplus 0) = \ind_{\VV}(L,\gamma_!(H\oplus 0)) = \ind(L,\CDS). \]
Taking these two equalities together, we obtain the statement of the theorem.
\endproof

\sub{Adding parameters.}
More generally, suppose that $\A$ also depends on a parameter $x$, that is,
$\A_x=(\Ami_x,\Ama_x)$ is a graph continuous family of Fredholm extension pairs.
Then the graph $\Gamma_x$ of $\Ami_x$, the graph $\Gamma'_x$ of $\Ama_x$
and the quotient $\VV_x = \Gamma'_x/\Gamma_x$ depend continuously on $x$.

If $L_x\subset\VV_x$ is a continuous family of closed subspaces (\q{boundary conditions}),
then the same reasoning as in the proof of Proposition \ref{prop:Gr-Reg-cont}
shows that the family $\A_L$ of realizations is graph continuous.

The Cauchy data space $\CDS_x$ of $\A_x$ also depends on $x$, 
but this dependence is not continuous in the general case.
In this paper we deal only with a specific case when the continuity of $\CDS$ is guaranteed.

\sub{Realizations of invertible extension pairs.}
Let us start with the simplest case when each pair $\A_x$ is \textit{invertible}
(that is, each $\Ami_x$ is injective and each $\Ama_x$ is surjective).
Then the subspace $H \oplus 0$ is transversal to $(\Gamma_x,\Gamma'_x)$ in $\Hh$.
By Proposition \ref{prop:tr-tr}, the Cauchy data space $\CDS = \gamma_!(H\oplus 0)$ 
depends continuously on $x$.

\begin{thm}\label{thm:ind-K0-inv}
Let $\A_x=(\Ami_x,\Ama_x)$ be a graph continuous family of invertible extension pairs
and $L=(L_x)$ be a continuous family of Fredholm boundary conditions for $\A_x$.
Then 
\[ \ind\A_L = \ind_{\VV}(L,\CDS)\in K^0(X). \] 
\end{thm}

\proof
The proof is the same as the proof of Theorem \ref{Thm:ind0}, 
but the reference to Theorem \ref{Thm:ind-Gr1} should be replaced 
by the reference to its family version, Theorem \ref{thm:ind-Gr}.
\endproof

\sub{Realizations of Fredholm extension pairs.}
More generally, let $\A_x=(A_x,A'_x)$ be a graph continuous family of \textit{Fredholm} extension pairs.

\begin{prop}\label{prop:C-cont}
	If the dimensions of $\ker \Ami_x$ and $\coker \Ama_x$ are locally constant,
	then the Cauchy data space $\CDS_x$ of $\A_x$ depends continuously on $x$. 
\end{prop}

\proof 
The case $k=k'=0$ is already explained in the previous subsection.
For arbitrary fixed $k$ and $k'$, the subspaces 
\[ K_x=\ker \Ami_x \txt{and} K'_x = H\ominus\ran \Ama_x \cong \coker \Ama_x \] 
depend continuously on $x$ and thus are fibers of the (finite-dimensional) 
subbundles $K$ and $K'$ of the trivial Hilbert bundle $H_X$, 
The orthogonal complements 
\[ \H = H_X\ominus K \txt{and} \H'=H_X\ominus K' = \ran \Ama \]
are locally trivial Hilbert bundles.
Let $\bar{\A} = (\bar{A},\,\bar{A}')$ be the restriction of $\A$ to the Hilbert bundle $\H$ 
considered as a pair of operators $\bar{A},\,\bar{A}'$ acting from $\H$ to $\H'$.
As we saw above, the Cauchy data space $\CDS_{\bar{\A}}$ of an invertible extension pair $\bar{\A}$ 
depends continuously on parameter.
Since $\VV_{\A}$ and $\VV_{\bar{\A}}$ are naturally isomorphic 
and the corresponding Cauchy data spaces coincide, we conclude that 
$\CDS_{\A}=\CDS_{\bar{\A}}$ depends continuously on parameter.
\endproof

\sub{Proof of Theorem \ref{Thm:ind-AL}.}
We keep notations of the proof of Proposition \ref{prop:C-cont}.
The realization $\A_L$ decomposes into the orthogonal sum 
$\A_L = \bar{\A}_L\oplus \bar{0}$, where $\bar{0}$ is the zero bundle homomorphism $K\to K'$.	
By the additivity of the index, $\ind\A_L = \ind(\bar{\A}_{L})+\ind(\bar{0})$.
Since the pair $\bar{\A}$ is invertible, 
Theorem \ref{thm:ind-K0-inv} implies $\ind(\bar{\A}_{L}) = \ind_{\VV}(L,\CDS)$.
Since $\ind\bar{0} = [K]-[K']$, we get $\ind\A_L = \ind_{\VV}(L,\CDS)+[K]-[K']$.
\endproof

Equivalently, Theorem \ref{Thm:ind-AL} can be deduced from Theorem \ref{thm:ind-Gr-KK'} 
in the same manner as Theorem \ref{thm:ind-K0-inv} is deduced from Theorem \ref{thm:ind-Gr}.

\section{The case of compact base space}\label{sec:alt}

This section is devoted to the case of a compact base space $X$.
In this case our theorems concerning $K^0$ index can be proved in a more standard way, as we will show now.
We present here proofs of Theorems \ref{Thm:ind-AL} and \ref{thm:ind-K0-inv}.
Theorems \ref{Thm:ind-Gr1} and \ref{thm:ind-Gr-KK'} are proven in exactly the same manner;
we omit their proofs in order to avoid repetition.

\sub{Stabilization and homotopy.}
Without restriction of generality, we can suppose in all our theorems
that the subspaces $L_x$ and $\VV_x/L_x$ are infinite-dimensional.
Indeed, if this is not the case, then one can use additivity of the index and stabilization 
in order to pass to the infinite-dimensional case, as described in Section \ref{sec:GrF}.

Let $\VV_x$, $\xX$ be a continuous family of subspaces of $\Hh$
and $(L_x,\CDS_x)$ be a continuous family of Fredholm pairs in $\VV_x$ 
such that $\dim \CDS_x=\codim \CDS_x=\infty$.
Since $X$ is compact, the index of the Fredholm pair $(L,\CDS)$ 
can be represented by the virtual difference of two vector bundles $V$, $W$ over $X$, 
\begin{equation}\label{eq:ind-LC}
	\ind_{\VV}(L,\CDS) = [V]-[W] \in K^0(X). 
\end{equation}
Embedding $V$ to $\CDS$ and $W$ to $\CDS^{\bot}=\VV\ominus \CDS$ as subbundles
and denoting $V'=\CDS\ominus V$, $W'=\CDS^{\bot}\ominus W$,
we define a Hilbert subbundle of $\VV$ 
\begin{equation*}%\label{eq:ind-L'}
 L'=V\oplus W' \subset \CDS\oplus \CDS^{\bot} = \VV. 
\end{equation*}
Since $L'$ is a finite rank deformation of $\CDS^{\bot}$, the pair $(L',\CDS)$ is Fredholm.
By the additivity of the index,
\[ \ind_{\VV}(L',\CDS) = \ind_{V\oplus W}(V,V) + \ind_{V'\oplus W'}(W',V') = [V]-[W]
   = \ind_{\VV}(L,\CDS)\in K^0(X), \]
since the pair $(W',V')$ is transversal in $V'\oplus W'$.
The equality $\ind_{\VV}(L',\CDS) = \ind_{\VV}(L,\CDS)$ implies that $L'$ is homotopic to $L$ 
(in the class of subbundles of $\VV$ that are Fredholm with respect to $\CDS$).

All our proofs below are based on this construction, and we will keep its notations.

\sub{Proof of Theorem \ref{thm:ind-K0-inv} for compact $X$.}
By the homotopy invariance of index, $\ind\A_{L} = \ind\A_{L'}$.
Since each $A_x$ is injective and each $A'_x$ is surjective, we have 
\[ \ker \A_{L'} = L'\cap\CDS = V \;\text{ and }\; \coker \A_{L'} \cong \VV/(L'+\CDS) = W, \;\text{ so}  \]
\[ \ind\A_{L} = \ind\A_{L'} = [\ker A_{L'}]-[\coker A_{L'}] = [V]-[W] = \ind_{\VV}(L,\CDS). \;\;\square\]

\sub{Proof of Theorem \ref{Thm:ind-AL} for compact $X$.}
By the assumption of the theorem, $\ker A_x$ and $\coker A'_x$ have locally constant dimension,
so they depend continuously on $x$ and give rise to locally trivial vector bundles 
$\ker A$ and $\coker A'$ over $X$.
Canonical exact sequences \eqref{eq:ex-sec1} and \eqref{eq:ex-sec2} applied to $L'$ instead of $L$ take the form
\begin{align*}
 &	0 \to \ker A \to \ker\A_{L'} \to L'\cap\CDS \to 0, \\ 
 &	0 \to \coker A' \to \coker\A_{L'} \to \VV/(L'+\CDS) \to 0. 
\end{align*}
Since $\ker A$ and $L'\cap\CDS=V$ are vector bundles, $\ker\A_{L'}$ is also a vector bundle and 
\[ [\ker\A_{L'}] = [\ker A]+[V]\in K^0(X). \]
Similarly, $\coker\A_{L'}$ is a vector bundle and \[ [\coker\A_{L'}] = [\coker A']+[W]. \]
Taking these two equalities together with \eqref{eq:ind-LC}, we obtain
\[ \ind\A_L = \ind\A_{L'} = [\ker\A_{L'}]-[\coker\A_{L'}] = \ind_{\VV}(L,\CDS)+[\ker A]-[\coker A'].  \;\;\square \]

\addtocontents{toc}{\vspace{\otstup}}

\section{Self-adjoint operators and Lagrangian Grassmanian}\label{sec:Lagr-sa}

\sub{Lagrangian Grassmanian.}
Let $\Lambda$ be a Hilbert space equipped with a skew-adjoint unitary operator $J$.
Then 
\begin{equation}\label{eq:J-omega}
	\omega(\xi,\eta) = \bra{J\xi,\eta}  
\end{equation}
is a symplectic form which turns $\Lambda$ into a (complex) symplectic space.
This means that $\omega$ is a non-degenerate sesquilinear form 
and $\omega(\eta,\xi) = -\overline{\omega(\xi,\eta)}$.

Recall that the annihilator (with respect to $\omega$) of a subspace $L\subset \Lambda$ is defined as
\[ L\om = \sett{\xi\in \Lambda}{\omega(\xi,\eta)=0 \text{ for every } \eta\in L}; \]
$L$ is called isotropic if $L\subset L\om$ and Lagrangian if $L=L\om$.
The part of the Grassmanian $\Gr(\Lambda)$ consisting of Lagrangian subspaces is called the Lagrangian Grassmanian of $\Lambda$ and denoted by $\LGr(\Lambda)$.

\sub{Unitary picture of the Lagrangian Grassmanian.}
The symplectic space $\Lambda$ is decomposed into the direct sum 
$\Lambda = \Lambda^+\oplus \Lambda^-$ of eigenspaces
\[ \Lambda^+ = \ker\br{J-i} \txt{and} \Lambda^- = \ker\br{J+i}. \]
As is well known, Lagrangian subspaces of $\Lambda$ are in one-to-one correspondence 
with unitary operators $u\colon \Lambda^+\to \Lambda^-$
(so $\Lambda$ admits Lagrangian subspaces if and only if $\dim\Lambda^+=\dim\Lambda^-$).
This correspondence takes a unitary operator $u$ to its graph 
$\gr_u\subset \Lambda^+\oplus \Lambda^- = \Lambda$ and provides a canonical homeomorphism 
\begin{equation}\label{eq:LGr-U}	
	\LGr(\Lambda)\isor\U(\Lambda^+,\Lambda^-). 
\end{equation}

\upskip
\sub{Self-adjoint operators and their graphs.}
Let $\Hhat=H\oplus H$ be the standard complex symplectic space with the symplectic form 
\begin{equation}\label{eq:omegaHH}
	\omega(\xi,\eta) = \bra{J\xi,\eta}, \quad\text{where}\; J=\smatr{0 & -1 \\ 1 & 0}. 
\end{equation}
Recall that $\Reg(H)$ is the space of regular operators equipped with the graph topology.
Let $\Rsa(H)$ denote its subspace consisting of self-adjoint operators.

The graph of the adjoint $A^*$ of a regular operator $A$ is the annihilator of the graph of $A$.
Therefore, $A$ is self-adjoint if and only if its graph $\gr_A$ is a Lagrangian subspace of $\Hh$.
This determines the canonical embedding 
\begin{equation*}%\label{eq:Rsa-LGr}
 \Rsa(H)\hookto\LGr(\Hh), \quad A\mapsto\gr_A, 
\end{equation*}
so $\Rsa(H)$ can be considered as a subspace of $\LGr(\Hh)$.

\sub{Cayley transform.}
Another well known embedding of $\Rsa(H)$, to the unitary group, is given by the Cayley transform 
\begin{equation}\label{eq:kap}
 \kap\colon\Rsa(H)\hookto\U(H), \quad \kap(A) = (A-i)(A+i)\inv. 
\end{equation}

The subspaces $\Hh^{\pm} = \ker(J\mp i)$ can be written explicitly as
\[ \Hh^+ = \sett{\xi\oplus i\xi}{\xi\in H} \txt{and} \Hh^- = \sett{\xi\oplus (-i\xi)}{\xi\in H}, \]
Identifying $\Hh^+$ and $\Hh^-$ with $H$ by the unitary operators 
\begin{equation}\label{eq:Hhat+-H}
	\xi\oplus (\pm i\xi) \mapsto \sqrt{2}\cdot\xi, 
\end{equation}
we obtain an identification of $\U(\Hh^+,\Hh^-)$ with $\U(H)$.
Denote by 
\[ \kab\colon\LGr(\Hh)\to\U(H) \] 
the composition of homeomorphism \eqref{eq:LGr-U} with this identification.

Let $L\in\LGr(\Hh)$ and $v=\kab(L)$. 
Then $L$ consists of vectors of the form
\[ (\xi\oplus i\xi) + (v\xi\oplus -iv\xi) = (1+v)\xi\oplus i(1-v)\xi \]
with $\xi$ running $H$.
Applying this to the graph $L=\gr_A$ of a regular self-adjoint operator $A$, we see that
$\dom A = \ran(1+v)$ and $i(1-v) = A(1+v)$, that is,
\[ \kab(\gr_A) = -(A+i)\inv(A-i) = -(A-i)(A+i)\inv = -\kap(A). \]
It follows that the Cayley transform can be extended continuously 
from the subspace $\Rsa(H)$ of $\LGr(\Hh)$ to the whole $\LGr(\Hh)$,
giving rise to the canonical homeomorphism
\begin{equation}\label{eq:kap-LGr}
  \kap\colon\LGr(\Hh)\to\U(H), \quad \kap(L)=-\kab(L). 
\end{equation}
(Such an extension is unique since $\Rsa(H)$ is dense in $\LGr(\Hh)$, but we don't need this fact.)
Taking all this together, we obtain the following commutative diagram:
\begin{equation}\label{diag:kap-kap}
	\begin{tikzcd}
	  \Rsa(H) \hookr{r} \hookr{rd}{\kap} 
		& \LGr(\Hh) \arrow{r} \arrow{d}{\kap} \arrow{rd}{\kab} & \U(\Hh^+,\Hh^-) \arrow{d}
		\\
		& \U(H) \arrow{r}[swap]{\times(-1)} & \U(H) 
	\end{tikzcd}
\end{equation}

\section{$K^1$ index}\label{sec:K1ind}

\upskip
\sub{Reduced unitary group.}
Let $H$ be an infinite-dimensional Hilbert space.
The unitary group $\U(H)$ equipped with the usual norm topology is a topological group,
with the multiplication defined by composition. 
The reduced unitary group 
\[ \UK(H) = \sett{u\in\U(H)}{u-1 \text{ is compact}} \] 
is the subgroup  of $\U(H)$ consisting of unitary operators which are compact deformation of the identity.
This group is well known to be a classifying space for the functor $K^1$.

The group $K^1(X)$ might be \textit{defined} as the set $[X,\UK]$ of homotopy classes of maps $X\to\UK(H)$.
The group operation in $K^1(X)$ is given by the group operation in $\UK$.

The subgroup $\UK(H)$ of $\U(H)$ is normal, that is, it is invariant under the conjugation by unitary operators.
The Kuiper's theorem implies that the class of $f\colon X\to\UK$ is invariant under the conjugation, that is, 
$[vfv\inv] = [f]$ for every continuous family of unitary operators $v\colon X\to\U(H)$.
Hence it does not matter which Hilbert space $H$ is chosen for the definition of $K^0(X)$:
one can use an arbitrary unitary operator $H'\to H$ for the identification of $[X,\UK(H')]$ with $[X,\UK(H)]$.

An equivalent way to define the group operation in $K^1(X)$, which will be more convenient to us,
uses the direct sum: 
\begin{equation}\label{eq:+K1}
	[f]+[g] = [f\oplus g] = [g\oplus f] \in \brr{X,\UK(H\oplus H)} = K^1(X) \quad\text{for } f,g\colon X\to\UK(H). 
\end{equation}

\upskip
\sub{The space $\UF$.}
The reduced unitary group is a subspace of 
\[ \UF(H) \;=\; \sett{u\in\U(H)}{u+1 \text{ is Fredholm}} \;\subset\; \U(H). \] 
Equivalently, $\UF(H)$ consists of unitary operators $u\in\U(H)$ such that $-1$ 
does not belong to the essential spectrum of $u$.

While this larger space $\UF(H)$ is not a group, it is still a classifying space for the functor $K^1$;
the subspace $\UK(H)$ is its deformation retract. 
This can be proven by the methods of \cite{ASi}, see the proof of \cite[Proposition 3.2]{ASi},
and was done by Kirk and Lesch in \cite[Lemma 6.1(1)]{KL}.
More precisely, it was shown by Kirk and Lesch that the embedding $\UK(H)\hookto\UF(H)$ 
is a weak homotopy equivalence.
A simple modification of their argument shows that $\UK$ is a deformation retract of $\UF$.
Indeed, the quotient map $\UF\to\UF/\UK$ has a section, 
so it is a trivial bundle with the structure group $\UK$.
The base $\UF/\UK$ of this bundle is contractible, so the fiber $\UK$ is a deformation retract of $\UF$.
An alternative proof of the homotopy equivalence of the embedding $\UK\hookto\UF$ is given in \cite{Pr21},
see \cite[Theorem C]{Pr21}.

It follows that a map $f\colon X\to\UF(H)$ determines the class in $K^1(X)$ which we call the index of $f$,
\[ \ind (f) = [f]\in [X,\UF] = [X,\UK] = K^1(X). \]
The same argument as for $\UK$ shows that 
\begin{equation}\label{eq:vfv-K1}
	\ind(vfv\inv) = \ind(f) \quad\text{for every continuous map } v\colon X\to\U(H). 
\end{equation}
Since $\UF$ is not a group, one can no longer use the group structure of $\U(H)$ 
in order to define the addition in the Abelian group $[X,\UF]$.
However, \eqref{eq:+K1} provides the following equality 
that may be used as a definition of the addition in $[X,\UF]$:
\begin{equation}\label{eq:UF+}
	[f]+[g] = [f\oplus g] \in[X,\UF(H\oplus H)] = K^1(X)  \quad\text{for } f,g\colon X\to\UF(H). 
\end{equation}
As before, the neutral element of the group $[X,\UF]$ is represented by the constant map $f_x\equiv 1$  
and the inverse of $[f]$ is represented by the map $x\mapsto f_x\inv=f_x^*$.

\sub{Stabilization.}
In the previous subsections we considered an infinite-dimensional Hilbert space $H$.
If $H$ is finite-dimensional, then $\U_K(H)=\U_F(H)=\U(H)$ is no longer a classifying space for $K^1$.
In the general case, one can use stabilization in order to define the index
of a map $f\colon X\to \UF(H)$ for a Hilbert space $H$ of \textit{arbitrary} dimension.
Namely, one replaces $f$ by the direct sum $f\oplus 1 \colon X\to\U_F(H\oplus H')$
of $f$ and the constant map $X\to\U_F(H')$,
where $H'$ is an auxiliary infinite-dimensional Hilbert space and $1$ denotes the identity in $\U(H')$.
For infinite-dimensional $H$, this definition agrees with the one given above due to equality \eqref{eq:UF+}.

If $d=\dim H$ is finite, then $\UF(H)=\U(H)\cong\U(\C^d)$.
In this case the class of a map $f\colon X\to \U(H)$ is the usual $K^1$ class 
of the map $X\to \U(\C^d)$ in the topological $K$-theory.

\sub{Fredholm Lagrangian Grassmanian.}
Let $\Lambda$ be a Hilbert space with symplectic form \eqref{eq:J-omega}, as in the previous section.
The Fredholm Lagrangian Grassmanian is defined as the Fredholm part of the Lagrangian Grassmanian:
\[ \LGr_F^2(\Lambda) = \sett{(L,M)\in\LGr(\Lambda)^2}{(L,M) \text{ is a Fredholm pair}}. \] 
\[ \LGr_F(\Lambda;M) = \sett{L\in\LGr(\Lambda)}{(L,M) \text{ is a Fredholm pair}}, \quad M\in\LGr(\Lambda). \] 
Let $u_M$ denote the unitary operator corresponding to $M$ 
under the homeomorphism $\LGr(\Lambda)\to\U(\Lambda^+,\Lambda^-)$.
If a Lagrangian subspace $M$ is fixed, then the unitary operator $u_M\colon \Lambda^+\to \Lambda^-$ 
can be used for identification of $\Lambda^+$ with $\Lambda^-$,
which provides a homeomorphism 
\[ \kap_M\colon\LGr(\Lambda)\to\U(\Lambda^+) \]
taking $M$ to $1$.
A pair $(L,M)$ of Lagrangian subspaces determines the unitary operator 
\[ \kap(L;M) := \kap_M(L) = u_M\inv u_L \in\U(\Lambda^+). \]
The pair $(L,M)$ is Fredholm if and only if the difference $u_L-u_M$ is Fredholm,
or, equivalently, $1-\kap(L;M)$ is Fredholm.
Thus we obtain a canonical homeomorphism
\begin{equation}\label{eq:LGrF-UF}
	\LGr_F(\Lambda;M)\to \UF(\Lambda^+), \quad L\mapsto -\kap(L;M) = -u_M\inv u_L = \kap(L;M^{\bot}) 
\end{equation}
and the induced homomorphism 
\begin{equation*}%\label{eq:LGrF-K1}
 \ind\colon [X,\LGr_F(\Lambda;M)]\to K^1(X), 
\end{equation*}
which is an isomorphism if $\dim\Lambda^+=\dim\Lambda^-=\infty$.

For the space $\LGr_F^2(\Lambda)$ of Fredholm Lagrangian pairs, 
there is a canonical homeomorphism
\begin{equation*}%\label{eq:LGrF2-UFU}
 \LGr_F^2(\Lambda) \to \UF(\Lambda^+)\times\U(\Lambda^+,\Lambda^-), 
   \quad (L,M)\mapsto (\kap(L;M^{\bot}),u_M) = (-u_M\inv u_L,\, u_M). 
\end{equation*}
The index of a family $(L,M)$ of Fredholm Lagrangian pairs $(L_x,M_x)$ is defined 
as the index of the first component of this map, that is,
\begin{equation}\label{eq:ind-LM}
 \ind(L,M) := \ind(\kap(L;M^{\bot})) = \ind(-\kap(L;M)) = \ind(-u_M\inv u_L) \in K^1(X). 
\end{equation}
If $\dim\Lambda^+ = \dim\Lambda^- = \infty$, then $\U(\Lambda^+,\Lambda^-)$ is contractible,
$[X,\UF(\Lambda^+)]\to K^1(X)$ is an isomorphism,
and the homotopy type of the family $(L,M)$ is defined by its index.

\sub{Standard symplectic space.}
Let $\Hhat=H\oplus H$ be the standard symplectic space with symplectic form \eqref{eq:omegaHH}.
The homeomorphism $\kap\colon\LGr(\Hh)\to\U(H)$
defined by formula \eqref{eq:kap-LGr} takes the horizontal subspace $H\oplus 0$ to $-1\in\U(H)$
and the vertical subspace $0\oplus H$ to $1\in\U(H)$.
Under the identification \eqref{eq:Hhat+-H} of $\Hh^+$ with $H$,
for every Lagrangian subspace $L$ we have the equality $\kap(L) = \kap(L;0\oplus H)$.
Therefore, for a family $L_x\in\LGr_F(\Hh;H\oplus 0)$, $\xX$ we get $\kap(L_x)\in\UF(H)$ and 
\begin{equation}\label{eq:indL-kap}
	\ind_{\Hh}(L,H\oplus 0) = \ind\kap(L) \in K^1(X).
\end{equation}

\upskip
\sub{Fredholm self-adjoint operators.}
By a classical result of Atiyah and Singer \cite{ASi}, 
the space $\BF\sa(H)$ of bounded Fredholm self-adjoint operators has three connected components.
Two of them, the space $\BFp(H)$ of essentially positive operators
and the space $\BFm(H)$ of essentially negative operators, are contractible, 
while the third component $\BF^{\star}(H)$ is a classifying space for the functor $K^1$.

In contrast with the bounded case, the space $\Rsa_F(H)$ of regular Fredholm self-adjoint operators 
equipped with the graph topology is path connected \cite[Theorem 1.10]{BLP}. 
As was shown by Joachim, this space is a classifying space for the functor $K^1$.
See \cite[Theorem 3.5(ii)]{Jo}.
An alternative proof of this result given by the author shows that both embeddings in the sequence 
\begin{equation*}%\label{eq:BF*-RFsa-LGrF}
  \BF^{\star}(H) \hookto \Rsa_F(H) \hookto \UF(H). 
\end{equation*}
are homotopy equivalences.
See \cite[Theorem C]{Pr21}.
Here the second embedding $\Rsa_F\hookto\UF$ is given by the Cayley transform \eqref{eq:kap}.
Therefore, one can \textit{define} the index of a graph continuous family $\Op = (\Op_x)_{\xX}$ 
of Fredholm self-adjoint operators parametrized by points of a topological space $X$
as the index of the composition $\kap\circ\Op\colon X\to\UF(H)$.
For bounded operators and compact base spaces this definition of the index 
is equivalent to the classical definition of Atiyah and Singer \cite{ASi}.

Taking into account equality \eqref{eq:indL-kap}, we obtain a $K^1$-analogue of formula \eqref{eq:indA2}:
\begin{equation}\label{eq:ind-op-gr}
	\ind\Op = \ind\kap(\Op) = \ind\kap(\gr_{\Op}) = \ind_{\Hh}(\gr_{\Op},H\oplus 0).
\end{equation}

\upskip
\sub{Properties.} 
We will use the following properties of the index in our proofs.

\begin{enumerate}[topsep=0pt, itemsep=0pt, parsep=3pt, partopsep=0pt]
	\item Homotopy invariance: $\ind_{\Lambda}(L,M)=\ind_{\Lambda}(L',M')$ 
	      for homotopic families of Fredholm Lagrangian pairs $(L_x,M_x)$ and $(L'_x,M'_x)$ in $\Lambda$.
	\item Additivity with respect to orthogonal sums: 
	      \[ \ind_{\Lambda\oplus\Lambda'}(L\oplus L',M\oplus M') = \ind_{\Lambda}(L,M)+\ind_{\Lambda'}(L',M'), \]
				for a family $(L_x,M_x)$ of Fredholm Lagrangian pairs in $\Lambda$ and $(L'_x,M'_x)$ in $\Lambda'$.
	\item Invariance: $\ind_{\Lambda}(L,M)=\ind_{\Lambda}(vL,vM)$ 
	      for a continuous family of unitary symplectomorphisms $v\colon X\to\U(\Lambda)$.
	\item Vanishing on families of transversal pairs.
\end{enumerate}

\medskip
The first two properties hold by the definition.
The third property follows from \eqref{eq:vfv-K1}. 
Indeed, $v$ commutes with $J$, so it decomposes into the direct sum 
$v=v_+\oplus v_-$ with respect to the decomposition $\Lambda=\Lambda^+\oplus \Lambda^-$.
Hence $\kap(vL,vM) = v_+\kap(L,M)v_+\inv$ and 
\[ \ind(vL,vM) = \ind(-\kap(vL,vM)) = \ind(-v_+\kap(L,M)v_+\inv) = \ind(-\kap(L,M)) = \ind(L,M).\]
The last property follows from the fact that a transversal Lagrangian pair $(L,M)$ 
is canonically homotopic (in the class of such pairs) to the pair $(M^{\bot},M)$.
Since $\kap(M^{\bot};M^{\bot})=1$, the index of such a family vanishes.
One can construct a canonical homotopy between $(L,M)$ and $(M^{\bot},M)$
similarly to the Fredholm Grassmanian case discussed in Section \ref{sec:K0ind}.
Indeed, a pair $(L,M)$ is transversal if and only if $L$ is the graph of a bounded operator 
$T_L\colon M^{\bot}\to M$.
If $L$ and $M$ are Lagrangian, then $T_L$ is self-adjoint.
A linear deformation of $T_L$ to zero operator leads to the deformation of $(L,M)$ to $(M^{\bot},M)$.

Alternatively, one can use functional calculus and the unitary picture of the Lagrangian Grassmanian.
A pair $(L,M)$ is transversal if and only if $1-\kap(L,M) = 1+\kap(L,M^{\bot})$ is invertible.
Let $\U_0$ denote the subspace of $\UF(\Lambda^+)$ consisting of unitary operators $u$ 
such that $1+u$ is invertible.
The functional calculus provides a contraction of this space to $1\in\U_0$, so $\ind(-\kap(L,M))=0$.

\sub{Symplectic spaces depending on parameter.}
Suppose that $J=J_x$ depends continuously on a parameter $\xX$.
Then 
\[ \Lambda_x^+ = \ker\br{J_x-i} = \ker\frac{iJ_x+1}{2} \]
is the kernel of the orthogonal projection $\frac{iJ_x+1}{2}$, 
so $\Lambda_x^+$ depends continuously on $x$.
The same holds for $\Lambda_x^- = \ker\frac{iJ_x-1}{2}$.

The symplectic space $\Lambda=\Lambda_x$ itself may also depend on $x$, say, 
be a fiber of a locally trivial symplectic Hilbert bundle. 
In this paper we will deal with a more specific situation, namely, 
with a fixed symplectic Hilbert space $(\Hhat,J)$
and a continuous family $\Lambda_x$ of $J$-invariant closed subspaces of $\Hhat$.
In this case, $\Lambda_x$ is also symplectic and 
the restriction $J_x$ of $J$ to $\Lambda_x$ depends continuously on $x$.
Since $\Lambda^{\pm}_x$ is a subspace of $\Hh^{\pm}$, 
the orthogonal projection onto $\Lambda^{\pm}_x = \Lambda_x\cap\Hhat^{\pm}$ is equal to $Q_xP^{\pm}$,
where $Q_x$ and $P^{\pm}$ denote the orthogonal projections of $\Hhat$ 
onto $\Lambda_x$ and onto $\Hhat^{\pm} = \ker\br{J\mp i}$ respectively. 
Therefore, $\Lambda^+_x$ and $\Lambda^-_x$ depend continuously on $x$.

Let $\Lambda^+$ and $\Lambda^-$ be the Hilbert bundles over $X$ 
with the fibers $\Lambda^+_x$ and $\Lambda^-_x$. 
The Hilbert bundle $\Lambda^+$ is a pullback of the canonical bundle over $\Gr(\Hhat^+)$ 
by the map $f\colon X\to\Gr(\Hhat^+)$, $x\mapsto\Lambda^+_x$.
Passing to a stabilization if necessary, we can suppose that the spaces 
$\Lambda^+_x$ and $\Hh/\Lambda^+_x$ are infinite-dimensional, so $f$ takes $X$ to $\Grinf(\Hhat^+)$.
As was explained in Section \ref{sec:Gr}, the canonical bundle over $\Grinf$ is trivial.
Therefore, $\Lambda^+$ is also trivial (at least after stabilization). 
The same holds for $\Lambda^-$. 

A continuous family of pairs $(L_x,M_x)$ of Lagrangian subspaces of $\Lambda_x$
gives rise a continuous family of unitary operators $\kap(L_x;M_x)\in\U(\Lambda^+_x)$.
Passing to a trivialization $\Lambda^+\cong X\times H$ takes $-\kap(L_x;M_x)$ 
to the continuous family of unitary operators $u_x\in\U(H)$ in the fixed Hilbert space $H$.
If the pairs $(L_x,M_x)$ are Fredholm, then $u_x\in\UF(H)$ and thus its index is defined.
We define the index $\ind(L,M)\in K^1(X)$ of the family $(L_x,M_x)$ as the index of the family $u_x$.
A different choice of a trivialization of $\Lambda^+$ leads to a conjugation of $u$ 
by a continuous family of unitary operators $v\colon X\to\U(H)$, which does not affect the index.
Therefore, our definition of $\ind(L,M)$ does not depend on the choice of a trivialization of $\Lambda^+$.

Clearly, all the properties of the index stated above for a fixed Hilbert space $\Lambda$ 
hold also for families $(L_x,M_x)\in\LGr_F^2(\Lambda_x)$ with $\Lambda_x$ depending on $x$.

\section{Index formula for Fredholm Lagrangian pairs}\label{sec:sympl}

Let $(\Hhat,\omega)$ be an infinite-dimensional symplectic space
with the symplectic form $\omega$ defined by a skew-adjoint unitary $J$ as in \eqref{eq:J-omega}.
Let $\Gamma$ be an isotropic subspace of $\Hhat$ and $\Gom = J\Gamma^{\bot}$ be its annihilator 
with respect to $\omega$.
We can now apply the discussion from Section \ref{sec:Gamma}
to the pair $\Gamma\subset\Gom$,
taking into account the additional symplectic structure on $\Hhat$.

\sub{Semi-Fredholm pairs.}
Let $M$ be a Lagrangian subspace of $\Hh$.
By \cite[Lemma IV.4.9]{Kato}, the sum $M\om+\Gom = M+\Gom$ is closed if and only if $M+\Gamma$ is closed.
Hence the following three conditions are equivalent:
\begin{enumerate}[topsep=-2pt, itemsep=0pt, parsep=3pt, partopsep=0pt]
	\item The pair $(M,\Gamma)$ is lower semi-Fredholm;
	\item The pair $(M,\Gom)$ is upper semi-Fredholm;
	\item $M$ is Fredholm with respect to $(\Gamma,\Gom)$.
\end{enumerate}

\sub{Quotient space.}
The quotient space $\VV = \Gom/\Gamma$ has the induced symplectic structure 
given by the formula 
\begin{equation}\label{eq:omega-V}
	\omega(\xi+\Gamma,\eta+\Gamma) = \omega(\xi,\eta). 
\end{equation}
As before, we will often identify $\VV$ with the orthogonal complement $\Gom\ominus\Gamma$.
Under this identification, the induced symplectic structure is given by the restriction of $\omega$ to $\VV$.

Both $\VV$ and $\Hh\ominus\VV = \Gamma\oplus J\Gamma$ are invariant under $J$, 
so the symplectic space $\Hhat$ is decomposed into the orthogonal sum of these two symplectic subspaces:
\[ \Hhat = \VV\oplus (\Gamma\oplus J\Gamma). \]

We will use the same notations as in Section \ref{sec:Gamma}:
$\gamma\colon\Gamma'\to\Gamma'/\Gamma=\VV$ is the quotient map
and $\gamma_!M = \gamma(M\cap \Gamma') \subset \VV$ for $M\subset\Hh$.

Clearly, a subspace $L\subset\VV$ is Lagrangian in $\VV$ if and only if
$\gammi(L) = L\oplus\Gamma$ is Lagrangian in $\Hhat$.

Conversely, suppose that $M\subset\Hh$ is Lagrangian.
This does not necessarily implies that $\gamma_!M$ is Lagrangian because $\gamma_!M$ is not necessarily closed.
However, if additionally $M+\Gamma$ is closed, then $\gamma_!M$ is Lagrangian.
Indeed, the annihilator of $\gammi (\gamma_! M) = M\cap\Gom+\Gamma$ is 
\[ (M\cap\Gom+\Gamma)\om = \overline{M+\Gamma}\cap\Gom = (M+\Gamma)\cap\Gom = M\cap\Gom+\Gamma \]
(here overline denotes the closure of a subspace),
so $\gammi (\gamma_! M)$ is Lagrangian.
It follows that $\gamma_!M$ is also Lagrangian.

\sub{Push-forward for transversal subspaces.} 
Consider the subspace
\[ \S = \sett{M\in\LGr(\Lambda)}{M \text{ is transversal to } (\Gamma,\Gom)} \]
of $\LGr(\Lambda)$.
For $M\in\S$ the sum $M+\Gamma$ is closed, so $\gamma_!$ takes $\S$ to $\LGr(\VV)$.
By Proposition \ref{prop:tr-tr}, the map $\gamma_!\colon\S\to\LGr(\VV)$ is continuous.

\begin{prop}\label{prop:SN}
  The fiber
	\[ \S_N = \sett{M\in\S}{\gamma_!M = N}, \] 
	of $\gamma_!\colon\S\to\LGr(\VV)$ over $N\in\LGr(\VV)$
	    is canonically homeomorphic to the Banach space 
	\[ \B\sa(\Gamma)\oplus\B(JN,\Gamma) \;\cong\; \B\sa(\Gamma)\oplus\B(N,\Gamma), \]
	with its zero element corresponding to $\bar{N}=N\oplus J\Gamma\in\S_N$.
\end{prop}

\proof
Recall that in the proof of Proposition \ref{prop:TN} we showed that
a subspace $M$ of $\Hhat$ is transversal to $(\Gamma,\Gom)$ and satisfies condition $\gamma_!M = N$
if and only if $M$ is the graph of a bounded operator $\mu\colon N\oplus F \to E\oplus\Gamma$ 
such that $\restr{P_E\mu}{N} = 0$.
Here $E=\VV\ominus N = JN$ and $F=\Hh\ominus\Gamma'=J\Gamma$,
so $\mu$ acts from $N\oplus J\Gamma$ to $JN\oplus\Gamma$.

Since $N\oplus J\Gamma$ and $JN\oplus\Gamma$ are orthogonal Lagrangian subspaces of $\Hhat$,
the graph of $\mu$ is Lagrangian if and only if the operator 
\[ \psi=\mu J \in\B(JN\oplus\Gamma) \] 
is self-adjoint.
The identity $P_{JN}\restr{\mu}{N} = 0$ takes the form $P_{JN}\restr{\psi}{JN} = 0$, so
\[ \psi = \matr{0 & g^* \\ g & f}, \txt{where} f\in\B\sa(\Gamma) \text{ and } g\in\B(JN,\Gamma). \]
Conversely, for each operator $\psi$ of this form, 
the graph of $\mu = \psi J\inv = -\psi J$ is an element of $\S_N$.
\endproof

\sub{Index formula for Fredholm Lagrangian pairs.}
The following result is a $K^1$ analogue of Theorem \ref{Thm:ind-Gr1}.

\begin{thm}\label{thm:ind-LGr1}
Let $M$ be a Lagrangian subspace and $\Gamma$ be a closed isotropic subspace of $\Hh$.
Suppose that $M$ is transversal to $(\Gamma,\Gom)$.
Let $L=(L_x)_{\xX}$ be a continuous family of Lagrangian subspaces of $\VV$ 
which are Fredholm with respect to $\gamma_!M$.
Then
  \[ \ind_{\Hh}(\gammi L,M) = \ind_{\VV}(L,\gamma_! M) \in K^1(X). \]
\end{thm}

\proof
The proof mostly repeats the proof of Theorem \ref{Thm:ind-Gr1}.
Let $N=\gamma_!M$.
By Proposition \ref{prop:SN}, $\S_N$ is path connected,
so there is a path $\M\colon[0,1]\to\S_N$ connecting $\M_0=M$ and $\M_1 = N\oplus J\Gamma$.
By the definition of $\S_N$, $\M_t$ is transversal to the pair $(\Gamma,\Gom)$ 
and $\gamma_!\M_t = N$ for every $t\in[0,1]$.
Since the pair $(L,\gamma_!\M_t) = (L,N)$ is Fredholm, 
the pair $(\gammi L,\M_t)$ is also Fredholm.
By the homotopy invariance of the index, 
\[ \ind(\gammi L,\,  M) = \ind(\gammi L,\,  \M_0) = \ind(\gammi L,\,  \M_1) 
    = \ind(\gammi L,\, N\oplus J\Gamma). \]
The Fredholm Lagrangian pair 
$(\gammi L,\, N\oplus J\Gamma) = (L\oplus \Gamma, \, N\oplus J\Gamma)$ in $\Hhat$
is the orthogonal sum of the Fredholm Lagrangian pairs $(L,N)$ in $\VV$
and $(\Gamma,J\Gamma)$ in $\Gamma\oplus J\Gamma$.
The second pair $(\Gamma,\, J\Gamma)$ is transversal, so its index vanishes.
Additivity of the index implies
\[ \ind(\gammi L,\, N\oplus J\Gamma) = \ind(L,\,N) + \ind(\Gamma,\,J\Gamma) = \ind(L,\,\gamma_! M)+0. \]
Taking these equalities together, we obtain
\[ \ind(\gammi L,\,  M) = \ind(\gammi L,\, N\oplus J\Gamma) = \ind(L,\,\gamma_! M). \quad\square\]

\sub{Adding parameters.}
More generally, suppose that $M$ and $\Gamma$ depend continuously on a parameter $\xX$.
Then $\VV_x=\Gom_x/\Gamma_x$ also depends continuously on $x$.
Identifying $\VV_x$ with $\Gom_x\ominus\Gamma_x$, we obtain the situation 
discussed in the end of the previous section.
Hence the index of a continuous family of Fredholm Lagrangian pairs $(L_x,N_x)$ in $\VV_x$ is well defined.

If $M_x$ is a continuous family of Lagrangian subspaces of $\Hh$ transversal to $(\Gamma_x,\Gom_x)$, 
then $\gamma_! M_x$ depends continuously on $x$ by Proposition \ref{prop:tr-tr}.
The same reasoning as in the proof of Theorem \ref{thm:ind-LGr1} now provides the following result.

\begin{thm}\label{thm:ind-LGr}
let $\Gamma_x,M_x\subset \Hhat$ and $L_x\subset\VV_x$
	be continuous families of closed subspaces parametrized by points of a topological space $X$.
	Suppose that $\Gamma$ is isotropic, $M$ is Lagrangian and transversal to $(\Gamma,\Gom)$,
	and $L$ is Fredholm with respect to $\gamma_!M$.
	Then
  \[ \ind_{\Hh}(\gammi L,M) = \ind_{\VV}(L,\gamma_! M) \in K^1(X). \]
\end{thm}

\section{Index formula for Fredholm self-adjoint extensions}\label{sec:K1-op}

Let $\Hhat=H\oplus H$ be the standard complex symplectic space with symplectic form \eqref{eq:omegaHH}.
Let $A$ be a symmetric operator on $H$ and $\Gamma$ be its graph.
The graph of $A^*$ is the annihilator $\Gom$ of $\Gamma$.
Symplectic form \eqref{eq:omega-A} on 
\begin{equation}\label{eq:K1-op-V}
	\VV=\dom A^* / \dom A \cong \Gom/\Gamma 
\end{equation}
coincides with the restriction of symplectic form \eqref{eq:omegaHH} to $\Gom\ominus\Gamma$.
The graph of $A_{L}$ coincides with $\gammi L$.
Therefore, discussion from the previous section can be applied to the graphs of operators.

In particular, $A$ is injective with closed range if and only if $A^*$ is surjective;
$A$ is lower semi-Fredholm if and only if $A^*$ is upper semi-Fredholm.
This allows to simplify our terminology.
Recall that an operator is called semi-Fredholm if it is either lower or upper semi-Fredholm.
A semi-Fredholm symmetric operator is necessarily lower semi-Fredholm, 
so we can omit the term \q{lower} in our assumptions for $A$.

Next, the Cauchy data space 
\begin{equation}\label{eq:K1-op-C}
 \CDS = \gamma(\ker A^*) = \gamma_!(H\oplus 0) 
\end{equation}
is a Lagrangian subspace of $\VV$ if the range of $A$ is closed.
The extension $A_L$ is Fredholm if and only if 
$A$ is semi-Fredholm and the pair $(L,\CDS)$ is Fredholm in $\VV$;
$A_L$ is invertible if and only if 
$A$ is injective, range of $A$ is closed, and the pair $(L,\CDS)$ is transversal in $\VV$.

We collect some of these properties in the following proposition.

\begin{prop}\label{prop:K1-op-L}
Let $A$ be a semi-Fredholm symmetric operator on $H$.
Then the Cauchy data space $\CDS = \gamma(\ker A^*)$  is a Lagrangian subspace of $\VV$
and Fredholm self-adjoint extensions of $A$ are parametrized 
by points of the Fredholm Lagrangian Grassmanian $\LGr_F(\VV;\CDS)$.
\end{prop}

Let $A_x$ be a graph continuous family of semi-Fredholm symmetric operators,
$\dim\ker A_x$ be locally constant, 
and $L_x$ be a continuous family of Fredholm self-adjoint boundary conditions for $A_x$.
Then, as we saw in Section \ref{sec:ind-ext}, 
$A_L$ is a graph continuous family of Fredholm self-adjoint operators.
Formula \eqref{eq:ind-op-gr} then takes the form
\begin{equation}\label{eq:ind-AL-gamma}
	\ind A_{L} = \ind_{\Hhat}(\Gamma_{A_L},\, H\oplus 0) 
	  = \ind_{\Hhat}(\gammi L,\, H\oplus 0) \in K^1(X). 
\end{equation}

\upskip
\sub{Proof of Theorem \ref{Thm:K1-AL}.}
Suppose first that $A$ is injective.
Then $H\oplus 0$ is transversal to $(\Gamma,\Gom)$ and Theorem \ref{thm:ind-LGr1} implies 
\[ \ind(\gammi L,H\oplus 0) = \ind(L,\gamma_!(H\oplus 0)) = \ind(L,\CDS). \]
Together with \eqref{eq:ind-AL-gamma} this implies the desired equality $\ind A_{L} = \ind(L,\CDS)$.

In the general case, consider the restriction $\bar{A}$ of $A$ to $H'=H\ominus K$, where $K=\ker A$.
Then $\VV_A$ and $\VV_{\bar{A}}$ are naturally isomorphic, the corresponding Cauchy data spaces coincide, 
and $A_L = \bar{A}_L\oplus 0$ with respect to the orthogonal decomposition $H=H'\oplus K$.	
By the additivity of the index, $\ind A_{L} = \ind\bar{A}_{L}$.
Since $\bar{A}$ is injective, the first part of the proof implies that
$\ind\bar{A}_{L} = \ind_{\VV}(L,\CDS)$ and thus $\ind A_{L} = \ind_{\VV}(L,\CDS)$. 
\endproof

\sub{Proof of Theorem \ref{Thm:K1-AL-k}.}
If all operators $A_x$ are injective, then $H\oplus 0$ is transversal to $(\Gamma_x,\Gom_x)$ for every $x$.
Applying Theorem \ref{thm:ind-LGr} to formula \eqref{eq:ind-AL-gamma}, we get
\[ \ind A_{L} = \ind(\gammi L,H\oplus 0) = \ind(L,\gamma_!(H\oplus 0)) = \ind(L,\CDS). \]
In the general case, consider the family of subspaces $K_x=\ker A_x$.
Since $\dim\ker A_x$ is locally constant, $K_x$ depends continuously on $x$.
Let $\bar{A}_x$ be the restriction of $A_x$ to $\H_x=H\ominus K_x$.
The subspaces $K_x$ and $\H_x$ give rise to the subbundles $K$ and $\H$ of the trivial Hilbert bundle $H_X$
(here $K$ has finite rank).
The family $\bar{A}$ acts on the fibers of $\H$ instead of a fixed Hilbert space.
Nevertheless,	our construction of the homotopy $\M$ is canonical,
so it works in fibers of the symplectic bundle $\H\oplus\H$ as well.
Thus we still obtain the equality $\ind(\bar{A}_{L}) = \ind_{\VV}(L,\CDS)$.
The $K^1$-index of zero operator $K\to K$ is equal to $0$.
Additivity of the index implies $\ind A_{L} = \ind \bar{A}_{L}$.
Taking these equalities together, we obtain $\ind A_L = \ind_{\VV}(L,\CDS)$.
\endproof

\section{Symplectic boundary value spaces}

In Sections \ref{sec:Lagr-sa}--\ref{sec:sympl} 
we were dealing with a symplectic structure \eqref{eq:J-omega}
given by a skew-adjoint \emph{unitary} operator $J$.
In particular, our definition of the index was given only in this framework.
However, it is often convenient to deal with boundary value spaces where 
Hilbert and symplectic structures are related in a more general manner.
In this section we describe the corresponding analogues 
of Theorems \ref{Thm:K1-AL} and \ref{Thm:K1-AL-k}.

\sub{Symplectic spaces.}
Let $\Lambda$ be a Hilbert space and $\Lambda'$ be its dual 
(that is, the space of all antilinear maps $\Lambda\to\C$).
A sesquilinear form $\omega\colon\Lambda\otimes\Lambda\to\C$ is called symplectic 
if $\omega(\eta,\xi) = -\overline{\omega(\xi,\eta)}$ and $\omega$ is a perfect pairing 
(that is, the map $\Lambda\to\Lambda'$, $\xi\to\omega(\xi,\cdot)$ is an isomorphism).
Such a pair $(\Lambda,\omega)$ is called a (complex) symplectic Hilbert space.
Equivalently, 
\[ \omega(\xi,\eta) = \bra{J\xi,\eta}, \] 
where $J$ is a skew-adjoint automorphism of $\Lambda$, not necessarily unitary.
If $J$ \textit{is} unitary, then we say that the inner product and the symplectic form on $\Lambda$ are \emph{compatible}. 

For a symplectic space $(\Lambda,\omega)$, one can always find a symplectomorphism 
(that is, a bijection preserving the symplectic form)
to a symplectic space $(\Lambda',\omega')$ such that $\omega'$ is compatible with the inner product on $\Lambda'$.
In fact, there is a canonical such isomorphism 
\begin{equation}\label{eq:SLambda}
	S\colon\Lambda\to\Lambda, \quad S=\abs{J}^{1/2}
\end{equation}
which takes $\omega$ to the symplectic form $\omega'$ compatible with the inner product.
Indeed,
\[ \omega'(\xi,\eta) = \omega(S\inv\xi,S\inv\eta) = \bra{S\inv JS\inv\xi,\eta} = \bra{J'\xi,\eta}, \] 
where $J' = S\inv JS\inv = J\abs{J}\inv$ is a skew-adjoint unitary commuting with $S$.
(Equivalently, one can replace the inner product on $\Lambda$ 
be the new inner product $\bra{\xi,\eta}' = \bra{S\xi,S\eta}$ compatible with $\omega$.)

\sub{Group of symplectomorphisms.}
Let $(\Lambda,\omega)$ be an infinite-dimensional symplectic space admitting Lagrangian subspaces.

\begin{prop}\label{prop:Sympl-contr}
The group $\Sympl(\Lambda)$ of symplectomorphisms of $\Lambda$ is contractible.
\end{prop}

\proof
Fix a transversal pair $(L,M)$ of Lagrangian subspace of $\Lambda$.
Applying a symplectomorphism $T$ to the pair $(L,M)$, we again obtain a transversal pair 
$L'=T(L)$, $M'=T(M)$ of Lagrangian subspaces.
Restriction of $\omega\colon\Lambda\otimes\Lambda\to\C$ to $L\otimes M$, resp. $L'\otimes M'$
provides a perfect pairing between $L$ and $M$, resp. $L'$ and $M'$.
Since $T$ preserves $\omega$, it takes the pairing $L\otimes M\to\C$ to $L'\otimes M'\to\C$.
Therefore, a symplectomorphism $T$ is defined by a transversal pair $(L',M')$ of Lagrangian subspaces
and an isomorphism $\restr{T}{L}\colon L\isor L'$,
while the second isomorphism $\restr{T}{M}\colon M\isor M'$ is determined by $\restr{T}{L}$.
Conversely, every such triple $(L',M',L\isor L')$ gives rise to a symplectomorphism of $\Lambda$.
It follows that $\Sympl(\Lambda)$ is a locally trivial fiber bundle 
over the space $\LGr_0(\Lambda)$ of transversal Lagrangian pairs $(L,M)$, 
whose fiber over $(L,M)$ is the subspace of $\B(L,L')$ consisting of invertible operators.
Since $\Lambda$ is infinite-dimensional, both the base space and fiber are contractible.
Therefore, the total space $\Sympl(\Lambda)$ is also contractible.
\endproof

\sub{Fredholm Lagrangian index.}
Let $(L_x,M_x)$, $\xX$ be a continuous family of Fredholm Lagrangian pairs in $(\Lambda,\omega)$.
We define its index as the index of the family $(SL,SM)$ in $(\Lambda,\omega')$, 
where $S$ is given by \eqref{eq:SLambda}.

The index depends only on the symplectic structure on $\Lambda$ and does not depend on its Hilbert structure, as the following proposition shows.

\begin{prop}\label{prop:K1ind-J}
	Let $T\colon (\Lambda,\omega)\to (\Lambda',\omega')$ be a symplectomorphism
	and $(L_x,M_x)$, $\xX$ be a continuous family of Fredholm Lagrangian pairs in $(\Lambda,\omega)$.
	Then $\ind_\Lambda(L,M) = \ind_{\Lambda'}(TL,TM)$.
	More generally, this equality holds also for a family $(T_x)$ of symplectomorphisms
	norm continuously depending on $x$.
\end{prop}

\proof
Without loss of generality, we can suppose that both $\omega$ and $\omega'$ are compatible with the corresponding inner products.
Using stabilization if necessary, we can suppose that both $\Lambda$ and $\Lambda'$ 
are infinite-dimensional.
Using decompositions $\Lambda=\ker(J-i)\oplus\ker(J+i)$ and $\Lambda'=\ker(J'-i)\oplus\ker(J'+i)$,
we see that there is an isometric symplectomorphism $(\Lambda',\omega')\to (\Lambda,\omega)$.
Composing it with $T$, we pass to the case $(\Lambda',\omega') = (\Lambda,\omega)$ 
and $T\in\Sympl(\Lambda)$.
In this case, Proposition \ref{prop:Sympl-contr} implies that $T$ is homotopic to the identity in $\Sympl(\Lambda)$,
and thus the family $(TL,TM)$ of Fredholm Lagrangian pairs in $\Lambda$ is homotopic to $(L,M)$
and their indices coincide.
\endproof

\sub{Symplectic space depending on parameter.}
Suppose now that the symplectic space itself depends norm continuously on $\xX$.
By that we mean the following: 
\begin{enumerate}[topsep=-2pt, itemsep=0pt, parsep=3pt, partopsep=0pt]
	\item the spaces $\Lambda_x$ are fibers of a locally trivial Hilbert bundle $\Lambda\to X$,
    whose structure group $\U$ is equipped with the norm topology, 
	\item the operators $J_x$ determined by the formula 
	$\omega_x(\xi,\eta) = \bra{J_x\xi,\eta}_{\Lambda_x}$ depend norm continuously on $x$.
\end{enumerate}
Then we say that $(\Lambda,\omega)$ is a symplectic bundle over $X$.

If this case, $S=\abs{J}^{1/2}\colon\Lambda\to\Lambda$ is a norm continuous bundle map.
If $(L_x,M_x)\in\LGr_F(\Lambda_x,\omega_x)$ is a norm continuous family of Fredholm Lagrangian pairs,
then we define its index as the index of the family $(S_xL_x,S_xM_x)$ in $(\Lambda,S_*\omega)$.

\begin{prop}\label{prop:K1ind-Jx}
Let $(\Lambda,\omega)$ and $(\Lambda',\omega')$ be symplectic bundles over a paracompact space $X$.
	Let $T_x\colon\Lambda_x\to\Lambda'_x$, $\xX$ be a norm continuous family of symplectomorphisms
	and $(L_x,M_x)\in\LGr(\Lambda_x,\omega_x)$ be a norm continuous family of Fredholm Lagrangian pairs.
	Then $\ind_\Lambda(L,M) = \ind_{\Lambda'}(TL,TM)$.
\end{prop}

\proof
Without loss of generality, we can suppose that the rank of $\Lambda$ and $\Lambda'$ is infinite 
and both $\omega$ and $\omega'$ are compatible with the corresponding inner products,
that is, the corresponding skew-adjoint bundle automorphisms 
$J\colon\Lambda\to\Lambda$ and $J'\colon\Lambda'\to\Lambda'$ are unitary.
Since both $J$ and $J'$ are norm continuous, $\ker(J\pm i)$ and $\ker(J'\pm i)$ 
are locally trivial Hilbert bundles.
Moreover, they have infinite rank and thus are trivial by Kuiper theorem.
Fixing their trivializations, we pass to the framework of Proposition \ref{prop:K1ind-J},
which provides the desired equality $\ind(L,M) = \ind(TL,TM)$.
\endproof

\sub{Another form of Theorems \ref{Thm:K1-AL} and \ref{Thm:K1-AL-k}.}
Let $A$ be a semi-Fredholm symmetric operator on $H$ and 
\begin{equation}\label{eq:dom2V}
	\Gamma\colon\dom A^*\to V 
\end{equation}
be a surjective map to a Hilbert space $V$ with $\ker\Gamma=\dom A$.
Then $\Gamma$ induces a natural isomorphism $T\colon\VV\to V$ 
which takes $\omega$ to a symplectic form $\omega_V$ on $V$, so that
\begin{equation}\label{eq:omega-AV}
	\bra{A^*\xi,\eta} - \bra{\xi,A^*\eta} = \omega_V(\Gamma\xi,\Gamma\eta)
		\quad\text{for } \xi,\eta\in\dom A^*. 
\end{equation}
Similarly, $T$ takes $\CDS = \gamma(\ker A^*)$ to the Lagrangian subspace $\cds = \Gamma(\ker A^*)$ of $V$.

In terms of the new \q{trace map} $\Gamma$, 
Fredholm self-adjoint extensions of $A$ are parame\-trized 
by the Fredholm Lagrangian Grassmanian $\LGr_F(V;\cds)$.
For $L\in\LGr_F(V;\cds)$ we denote by $A_L$ the restriction of $A^*$ to the domain
$\dom A_L = \Gamma\inv(L)$.
Using Proposition \ref{prop:K1ind-J}, we can now state Theorem \ref{Thm:K1-AL} as follows.

\begin{thm}\label{thm:K1-V}
Let $L=(L_x)_{\xX}$, $L_x\in\LGr_F(V;\cds)$ be a continuous family of Fredholm self-adjoint boundary conditions for $A$.
Then $\ind A_L = \ind_{V}(L,\cds)\in K^1(X)$.
\end{thm}

More generally, let $A=(A_x)$ be a graph continuous family of semi-Fredholm symmetric operators
parametrized by points of a paracompact space $X$. 
Then the domains $\dom A^*_x$ form fibers of a locally trivial Hilbert bundle with the structure group $\U$.
Let $V\to X$ be another locally trivial Hilbert bundle, 
and let $\Gamma_x\colon\dom A^*_x\to V_x$ be \q{trace maps} as above 
which depend norm continuously on $x$.

\begin{thm}\label{thm:K1-Vx}
Suppose that $\dim\ker A_x$ is locally constant.
Then $\cds_x = \Gamma(\ker A^*_x)\subset V_x$ depends continuously on $x$ and
$\ind A_L = \ind_{V}(L,\cds)\in K^1(X)$ for every continuous family $L=(L_x)$
of Fredholm self-adjoint boundary conditions $L_x\in\LGr_F(V_x;\cds_x)$.
\end{thm}

\proof
$\Gamma$ induces a norm continuous bundle isomorphism $T\colon\VV\to V$.
The symplectic structure $J_x$ on $V_x$ is induced by this isomorphism 
from the symplectic structure on $\VV_x$, so it depends norm continuously on $x$.
Since $\dim\ker A_x$ is locally constant, $\CDS_x\subset\VV_x$ and thus also $\cds_x=T_x(\CDS_x)$
depend continuously on $x$.
Combining Theorem \ref{Thm:K1-AL-k} and Proposition \ref{prop:K1ind-Jx}, we obtain
$\ind A_L = \ind_{\VV}(T\inv(L),\CDS) = \ind_V(L,\cds)$.
\endproof

\section{Boundary triplets}\label{sec:triplets}

\upskip
\sub{Boundary triplets.}
Let again $A$ be a semi-Fredholm symmetric operator on $H$.
A \emph{boundary triplet} for $A$ is a triple $(K,\Gamma_0,\Gamma_1)$,
where $K$ is a Hilbert space and $\Gamma_0,\Gamma_1\colon \dom A^*\to K$ 
are bounded operators such that
\[ \Gamma = \Gamma_0\oplus\Gamma_1\colon \dom A^*\to K\oplus K \]
is a surjective operator satisfying the Green formula in the following form:
\begin{equation}\label{eq:Green-triplet}
	\bra{A^*\xi,\eta}_H - \bra{\xi,A^*\eta}_H 
	  = \bra{\Gamma_1\xi,\Gamma_0\eta}_K - \bra{\Gamma_0\xi,\Gamma_1\eta}_K
		\txt{for} \xi,\eta\in\dom A^*. 
\end{equation}
In this case, the kernel of $\Gamma = \Gamma_0\oplus\Gamma_1$ coincides with $\dom A$, 
so $\Gamma$ induces an isomorphism 
\[ T\colon \VV = \dom A^*/\dom A \to K\oplus K \txt{such that} \Gamma=T\circ\gamma. \]
Moreover, \eqref{eq:Green-triplet} implies that $T$ takes a symplectic structure on $\VV$ to the 
standard symplectic structure \eqref{eq:omegaHH} on $K\oplus K$.
It follows that a boundary triplet is given by a (not necessarily isometric) symplectomorphism
\[ T=T_0\oplus T_1\colon\VV\to K\oplus K, \txt{with} \Gamma_i=T_i\circ\gamma. \]
In these terms, self-adjoint boundary conditions for $A$
are given by Lagrangian subspaces $L$ of $K\oplus K$, so that
\begin{equation}\label{eq:AL-triplet}
	\dom A_L = \sett{\xi\in\dom A^*}{\Gamma\xi\in L}.
\end{equation}
Theorems \ref{thm:K1-V} and \ref{thm:K1-Vx} then take the following form.

\begin{thm}\label{thm:K1-triplet}
Let $(K,\Gamma_0,\Gamma_1)$ be a boundary triplet for a regular symmetric operator $A$, 
$\cds = \Gamma(\ker A^*)$, and $L=(L_x)_{\xX}$ 
be a continuous family of self-adjoint Fredholm boundary conditions $L_x\subset K\oplus K$ for $A$.
Then $\ind A_L = \ind_{K\oplus K}(L,\cds)\in K^1(X)$.

More generally, the same equality holds if $A_x$ depends graph continuously on $x$,
$\dim A_x$ is locally constant, 
and the boundary triplet $(K_x,\Gamma_x)$ (and then also $\cds_x$) depends continuously on $x$.
\end{thm}

\sub{Example: quantum graphs}
Let $G$ be a finite metric graph whose edges $e$ are intervals of length $\ell_e$.
Let $D=-d^{2}/dt^{2}+q(t)$ be the Schr\"odinger operator on $G$ with a potential $q\colon G\to\R$.
The operator $D$ is formally self-adjoint and its natural maximal domain is the second order Sobolev space 
\[ \dom D\ma = \bigoplus_{e}H^{2}(e) \subset L^2(G). \]
Let $\dG = \bigoplus_{e}\partial e$ be the disjoint union of the ends of all edges and
\[\Gamma_0,\Gamma_1\colon \dom D\ma \to K=L^{2}(\dG) = \C^{\abs{\dG}},
\quad \Gamma_0 f=f|_{\dG}, \quad \Gamma_1 f = \nabla f|_{\dG} \]
be the Dirichlet and Neumann trace maps.
Here $\nabla f$ denotes the derivative of $f$ at the end of an edge taken in  
the direction coming out from the vertex into the edge.
Then the adjoint operator $D\mi = (D\ma)^*$ has the domain $\dom D\mi = \ker\Gamma$ 
and the Green formula \eqref{eq:Green-triplet} is satisfied, so 
$(K,\Gamma_0,\Gamma_1)$ is a boundary triplet for $D$.

We will state the index formula for the most general \q{global} boundary conditions,
which are given by subspaces $L$ of $K\oplus K$.
Since $K$ is finite-dimensional, every boundary condition $L\subset K\oplus K$ is Fredholm.
It is self-adjoint if $L$ is Lagrangian in $K\oplus K$.
Theorem \ref{thm:K1-triplet} applied to $D$ takes the following form:

\begin{thm}\label{thm:index-QG}
Let $L=(L_x)_{\xX}$ be a continuous family of Lagrangian subspaces $L_x\subset K\oplus K$.
Then $\ind D_L = \ind_{K\oplus K}(L)\in K^1(X)$. 
\end{thm}

\upskip
\sub{Remarks.}
1. Since $K$ is finite-dimensional, the $K^1$-index of a family $(L_x,\CDS)$ is independent of $\CDS$,
so we write $\ind(L)$ instead of $\ind(L,\CDS)$.

2. Boundary conditions for $D$ are usually local, that is, 
are given separately at each vertex $v$ of the graph 
by a Lagrangian subspace $L^v$ of $K^v\oplus K^v$, 
where $K^v$ is the vector space whose dimension is the degree of the vertex $v$.
Then $K=\oplus_v K^v$, $L=\oplus_v L^v$,
and the theorem takes the form $\ind D_L = \oplus_v\ind(L^v)$, 
where the index of $L_v$ is taken in $K^v\oplus K^v$. 

3. A particular case of Theorem \ref{thm:index-QG} is 
the \q{spectral flow = Maslov index} formula for a finite quantum graph.
Indeed, if $X$ is a circle, then, under the identification  $K^1(S^1)\cong\Z$, 
the index of $L$ is the Maslov index of this loop of Lagrangian subspaces
and the index of $D_L$ is the spectral flow of this loop of self-adjoint operators.
This particular case was proven in \cite{BPS} by different methods. 
See \cite[Theorem A]{BPS}.

\sub{Boundary triplets and transversal pairs of Lagrangian subspaces.}
A boundary triplet $T$ determines a transversal pair $(\ker T_1,\ker T_0)$
of Lagrangian subspaces of $\VV$.
Conversely, every transversal pair $(K,N)$ of Lagrangian subspaces of $\VV$ determines a boundary triplet 
\[ T = T_0\oplus T_1\colon\VV\to K\oplus K \txt{with} K=\ker T_1 \txt{and} N=\ker T_0. \]
Indeed, comparing the perfect pairing $\omega\colon K\otimes N\to\C$ 
and the inner product $K\otimes K\to\C$, we obtain a natural isomorphism 
$\varphi\colon N\to K$ such that 
\[ \omega(\xi,\eta) = \bra{\xi,\varphi\eta} \txt{for all} \xi\in K, \eta\in N. \]
The direct sum decomposition $\VV=K\dotplus N$ leads to the symplectic isomorphism 
\[ \VV = K\dotplus N  \;\stackrel{\Id\oplus\varphi}{\longrightarrow}\; K\oplus K. \]
Therefore, boundary triplets $(K,\Gamma)$ are, up to a (probably non-isometric) automorphism of $K$,
in one-to-one correspondence with transversal pairs of Lagrangian subspaces in $\VV$.
In particular, fixing a closed subspace $N$ of $\VV$ transversal to $\CDS$,
we obtain the boundary triplet
\begin{equation*}%\label{eq:CN-CC}
	T\colon \VV = \CDS\dotplus N \isor \CDS\oplus\CDS. 
\end{equation*}

\sub{Boundary triplets given by invertible extensions.}
Suppose now that a symmetric operator $A$ is injective.
Let again $N$ be a closed subspace of $\VV$ transversal to $\CDS$.
Then the extension $\At=A_N$ is invertible. 

The choice of an invertible extension $\At$, which is sometimes called a \emph{reference operator}, 
leads to the following construction of a boundary triplet
\cite[Example 14.6]{Schm}. 
First, there are direct sum decompositions
\begin{equation}\label{eq:triplet-dom}
	\dom A^* = \dom\At\dotplus K \txt{ and } \dom\At = \dom A \dotplus \At\inv K, \txt{where} K=\ker A^*.
\end{equation}
Taken together, they provide the decomposition
\[ \dom A^* = \dom A\dotplus V, \txt{where} V = K\dotplus \At\inv K \subset H. \]
The projection $\dom A^*\to V$ onto the second summand is a \q{trace map} for $A$, as in \eqref{eq:dom2V}.

Let $\xi_i,\eta_i\in K$, $i=0,1$. 
Then $\xi=\xi_0+\At\inv\xi_1, \; \eta=\eta_0+\At\inv\eta_1 \in V\subset\dom A^*$ and
\[ 	\bra{A^*\xi,\eta}_H - \bra{\xi,A^*\eta}_H 
	  = \bra{\xi_1,\eta_0+\At\inv\eta_1}_H - \bra{\xi_0+\At\inv\xi_1,\eta_1}_H 
		= \bra{\xi_1,\eta_0}_K - \bra{\xi_0,\eta_1}_K, 
\]
since $\At\inv$ is self-adjoint.
Therefore, the map 
\begin{equation}\label{eq:triplet-At}
	\Gamma = 0\dotplus\Id\dotplus\At \colon \dom A^* = \dom A\dotplus K\dotplus \At\inv K \to K\oplus K 
\end{equation}
provides the boundary triplet $(K,\Gamma)$ for $A$.

The trace map $\gamma\colon\dom A^*\to\VV$ takes $K$ to $\CDS$. 
Since $\At=A_N$ and thus $\dom\At = \gamma\inv(N)$, 
the second equality of \eqref{eq:triplet-dom} implies $\gamma(\At\inv K) = N$.
Therefore, the restriction of $\gamma$ to $V = K\dotplus \At\inv K$
takes the transversal Lagrangian pair $(K,\At\inv K)$ in $V$
to the transversal Lagrangian pair $(\CDS,N)$ in $\VV$,
and we obtain a commutative diagram of isomorphisms
\begin{equation*}%\label{diag:trip}
	\begin{tikzcd}
		V = K\dotplus \At\inv K \arrow{r}{\gamma} \arrow[swap]{d}{\Id\oplus\At}
		 & \CDS\dotplus N = \VV \arrow{ld}{T}  
		\\
		K\oplus K & 
	\end{tikzcd}
\end{equation*}
where $T(\CDS)=K\oplus 0$ and $T(N)=0\oplus K$.

As in \eqref{eq:AL-triplet}, for a Lagrangian subspace $L\subset K\oplus K$ 
we denote by $A_L$ the self-adjoint extension of $A$ with the domain 
$\dom A_L = \sett{\xi\in H}{\Gamma\xi\in L}$, where $\Gamma$ is given by \eqref{eq:triplet-At}.
Note that $A_L$ depends not only on $L$ but also \emph{on the choice of a reference operator} $\At$.

Lagrangian subspaces of $K\oplus K$ are often called \emph{self-adjoint relations} in $K$
and are considered as a generalization of self-adjoint operators $K\to K$.
Similarly, closed subspaces of $K\oplus K$ which are Fredholm with respect to $K\oplus 0$
are called \emph{Fredholm relations} in $K$.
The index $\ind L\in K^1(X)$ of a family $L_x$ of Fredholm self-adjoint relations 
is defined as the index of the family $(L, K\oplus 0)$ of Fredholm Lagrangian pairs in $K\oplus K$.
Equivalently, it is equal to the index of the Cayley transform $\kap(L)$, see \eqref{eq:indL-kap}.

Suppose now that both $A=A_x$ and $\At=\At_x$ depend graph continuously on $x$
and that each $\At_x$ is invertible, as above.
These data lead to a continuous family of boundary triplets $(K_x,\Gamma_x)$, where $K_x=\ker A^*_x$.
In this framework, Theorem \ref{thm:K1-triplet} takes the following form.

\begin{thm}\label{thm:K1-triplet-inv}
Let $A_x$ and $\At_x$ satisfy assumptions above
and $L_x$ be a continuous family of Fredholm self-adjoint relations in $K_x = \ker A^*_x$.
Then $\ind A_L = \ind L \in K^1(X)$.
\end{thm}


\begin{thebibliography}{9999}

\bibitem[A]{Atiyah}
M. F. Atiyah.
K-theory, Lecture notes by DW Anderson. 
Harvard, Benjamin, New York (Fall, 1964) (1967).

\bibitem[AS]{ASi}
M.F. Atiyah and I. M. Singer. 
Index theory for skew-adjoint Fredholm operators. 
Publications math\'ematiques de l'IH\'ES \textbf{37} (1969), no.1, 5--26.

\bibitem[BPS]{BPS} 
R. Band, M. Prokhorova, and G. Sofer.
Spectral flow and Robin domains on metric graphs.
To appear in the Israel Journal of Mathematics;
arXiv:2505.02039 (2025), 33 pp.

\bibitem[BGS]{BGS}
L. Bandara, M. Goffeng, and H. Saratchandran. 
Realisations of elliptic operators on compact manifolds with boundary. 
Advances in Mathematics \textbf{420} (2023): 108968.

\bibitem[BLP]{BLP}
 B. Booss-Bavnbek, M. Lesch, and J. Phillips.
 Unbounded Fredholm operators and spectral flow.
 Canadian Journal of Mathematics \textbf{57} (2005), no.2, 225--250;
 arXiv:math/0108014 [math.FA].

\bibitem[CL]{CL}
H. O. Cordes, J. P. Labrousse.
The invariance of the index in the metric space of closed operators.
Journal of Mathematics and Mechanics \textbf{12} (1963), no.5, 693--719.

\bibitem[DS]{DS}
N. Dunford and J.T. Schwartz. 
Linear operators II: spectral theory (1964).

\bibitem[$I_1$]{I23}
N.V. Ivanov. 
Boundary triplets and the index of families of self-adjoint elliptic boundary problems. 
arXiv:2306.15132 (2023), 45 pp.

\bibitem[$I_2$]{I23b}
N.V. Ivanov. 
On the index of families of self-adjoint abstract boundary problems.
Preprint, December 2023, 9 pp.

\bibitem[J]{Jo}
M. Joachim.
Unbounded Fredholm operators and K-theory.
High-dimensional Manifold Topology, World Sci. Publ., River Edge, NJ (2003), 177--199.

\bibitem[Ka]{Kato}
T. Kato. 
Perturbation Theory for Linear Operators. 
A Series of Comprehensive Studies in Mathematics \textbf{132} (1980).

\bibitem[KL]{KL}
P. Kirk, M. Lesch. 
The h-invariant, Maslov index, and spectral flow for Dirac-type operators on manifolds with boundary. 
Forum Math \textbf{16} (2004), 553--629.

\bibitem[Ku]{Kui}
N.H. Kuiper. 
The homotopy type of the unitary group of Hilbert space. 
Topology, \textbf{3} (1965), no. 1, 19--30.

\bibitem[N]{Neu}
G. Neubauer. 
Homotopy properties of semi-Fredholm operators in Banach spaces. 
Mathematische Annalen \textbf{176} (1968), no.4, 273--301.

\bibitem[$P_1$]{Pr17}
M. Prokhorova.
Self-adjoint local boundary problems on compact surfaces. I. Spectral flow.
Journal of Geometric Analysis, \textbf{31} (2021), no. 2, 1510--1554;
arXiv:1703.06105.

\bibitem[$P_2$]{Pr21}
M. Prokhorova.
Spaces of unbounded Fredholm operators. I. Homotopy equivalences.
Journal of Topology and Analysis, 18 (2026), no. 4, 1213-1241;
arXiv:2110.14359.

\bibitem[RS]{RS}
M. Reed and B. Simon. 
Methods of modern mathematical physics. II: Fourier analysis, self-adjointness. Elsevier, 1975.

\bibitem[S]{Schm}
K. Schm\"udgen. 
Unbounded self-adjoint operators on Hilbert space. 
Springer, 2012, 432 pp.

\bibitem[WO]{WO}
N.E. Wegge-Olsen.
K-theory and C*-algebras: a Friendly Approach. Vol. 1050. Oxford: Oxford university press, 1993.

\end{thebibliography}
\end{document}